\documentclass[12pt]{amsart}
\usepackage[utf8]{inputenc}
\usepackage{amsfonts, bm}
\usepackage{amsbsy}
\usepackage{amssymb}
\usepackage{amsxtra}
\usepackage{amsmath}
\usepackage{graphics}
\usepackage{enumerate}
\usepackage[usenames]{color}
\usepackage[english]{babel}
 \newcommand{\obra}[3]{{\sc #1} {\em #2}. {#3}.}

\newtheorem{theorem}{\bf Theorem}
 \newtheorem{lemma}[theorem]{\bf Lemma}
 \newtheorem{proposition}[theorem]{\bf Proposition}
 \newtheorem{definition}[theorem]{\bf Definition}
 
 \newtheorem{remark}[theorem]{Remark}

\renewenvironment{proof}{{\em Proof \/.-}}
   {\hfill $\square$\newline}
\newcommand{\R}{\mathbb{R}}
 
 \newcommand{\N}{\mathbb{N}}
 
  \newcommand{\C}{\mathbb{C}}

 \newcommand{\CC}{\mathcal{C}}

    \newcommand{\MM}{\mathcal{M}}

 \newcommand{\xx}{{\bm x}}
 \newcommand{\yy}{{\bm y}}
 \newcommand{\zz}{{\bm z}}
  \newcommand{\ww}{{\bm w}}

    \newcommand{\rr}{{\bm r}}
      \renewcommand{\ss}{{\bm s}}
  \newcommand{\vp}{\varphi}
  \newcommand{\wt}[1]{{\widetilde{#1}}}
   \newcommand{\wh}[1]{{\widehat{#1}}}

\newcommand{\g}{\gamma}
\newcommand{\G}{\Gamma}

\newcommand{\GL}{\text{GL}}

\newcommand{\ord}{\text{ord}}
\renewcommand{\wt}[1]{\widetilde{#1}}
\renewcommand{\wh}[1]{\widehat{#1}}

\begin{document}
	
	\title[Trajectories asymptotic to formal curves]{Trajectories of vector fields asymptotic to formal invariant
		curves}
	\author{O. LeGal}
	\address{Université Savoie Mont Blanc. Laboratoire de Mathématiques, LAMA. Bâtiment Le Chablais, Campus Scientifique, 73376 Le Bourget du Lac, France.}
	\email{Olivier.Le-Gal@univ-smb.fr}
	\author{F. Sanz S\'{a}nchez}
	\address{Universidad de Valladolid. Departamento de \'{A}lgebra, An\'{a}lisis Matem\'{a}tico, Geometr\'{\i}a y Topolog\'{\i}a. Facultad de Ciencias. Paseo de Belén, 7, E-47011, Valladolid (Spain)} \email{fsanz@uva.es}
	\thanks{Second author was partially supported by the Agencia Estatal de Investigación, Ministerio de Ciencia e Innovaci\'on (Spain). Project PID2019-105621GB-I00}

\begin{abstract} {We prove that a formal curve $\Gamma$ that is invariant by a $C^{\infty}$ vector field $\xi$ of $\mathbb R^{m}$ has a geometrical realization, as soon as the Taylor expansion of $\xi$ is not identically zero along $\Gamma$. This means that 
there is a trajectory $\gamma\subset \mathbb R^{m}$ of $\xi$ which is asymptotic to $\Gamma$. 
This result solves a natural question proposed by Bonckaert nearly forty years ago. 
We also construct an invariant $C^0$ manifold $S$ in some open horn around $\Gamma$ which is composed entirely of trajectories asymptotic to $\Gamma$, and contains the germ of any such trajectory. 
If $\xi$ is analytic, we prove that there exists a trajectory $\gamma$ asymptotic to $\G$ which is, moreover, non-oscillating with respect to subanalytic sets.}
\end{abstract}


\maketitle 

\subsection*{Acknowledgement}
This project is born during an enriching discussion with Felipe Cano. We would like to thank him for that and for his precious advises during the progress of this work.

\section{Introduction}

In this article, we consider germs of smooth vector fields $\xi \in \text{Der}_{\mathbb R}(C^{\infty}(\mathbb R^{m},0))$
which admit an irreducible formal invariant curve, this is, a $\Gamma\in\mathbb (R[[t]])^{m}$ such that $(\hat\xi \circ \Gamma) \wedge\Gamma'=0$, where $\hat\xi\in\text{Der}_{\mathbb R}(\mathbb R[[x_1,\dots,x_{m}]])$ is the Taylor expansion of $\xi$.
The ordinary correspondance between formal series and asymptotic development incites to believe there should be a a real curve $\gamma \subset\mathbb R^{m}$ invariant by $\xi$ and asymptotic to $\Gamma$. {The main result in this paper proves that this intuition is actually true, under the hypothesis (in general necessary) that $\Gamma$ is not included in the formal singular locus of $\xi$; that is, $\hat\xi \circ \Gamma\neq 0$.}

\begin{theorem}
\label{thm:1}
Let $\xi$ be the germ of a $C^{\infty}$ vector field at 
$a\in \R^{m}$, and  $\G$ a formal curve at $a$, invariant by $\xi$, 
that is not included in the formal singular locus of $\xi$. 
Then, there exists a germ of trajectory $\g$ of $\xi$ which is
asymptotic to $\Gamma$. 
\end{theorem}

In \cite{Bon}, Bonckaert solves the $3$ dimensional case and asks whether this very natural 
result can be generalized to higher dimension (see \cite[Remark 2.4, p. 118]{Bon}).
A broader question, already addressed in that paper in dimension three, is to describe {a realization of the ``attracting basin" of $\Gamma$}. We prove that such a realization exists as a topological manifold. In fact, there is one such manifold associated to each of the two half-branches of $\Gamma$ (the formal analogous of the connected components of $\G\setminus\{0\}$ when $\G$ converges). We postpone the precise definitions and a more accurate statement to section \ref{sec:final}:
 
\begin{theorem}\label{thm:2}
Let $\xi$ be a $C^{\infty}$ vector field at 
$a\in \R^{m}$, and $\G$ a formal curve at $a$, invariant by $\xi$, 
such that $\wh{\xi}\circ\G\ne 0$. 
{
For each half-branch $\G^+$ of $\G$ there is a finite order horn neighborhood $V^+$ of $\G^+$ 
and a connected topological submanifold $S^+\subset V^+$ of positive dimension with the following property. For any $b\in V^+$, the trajectory of $\xi$ issued from $b$ 
is asymptotic to $\G^+$ iff $b\in S^+$, and escapes from $V^+$ otherwise.}
\end{theorem}
 
Over the field of complex numbers $\mathbb C$, one usually considers holomorphic vector fields
and the theory of Borel-Laplace multi-summability answers the analogue question.
Indeed, a (complex) formal series $\Gamma$ that is invariant by a (complex) analytic vector field $\xi$,
and not contained in $\xi^{-1}(0)$, can be proven to be of Gevrey type and multi-summable. By a summation process 
as that proposed, among others, by Balser \cite{Bal}, Braaksma \cite{Bra}, Ramis \cite{Ram} or Malgrange \cite{Mal} 
we get invariant complex curves $\gamma$, {defined and asymptotic to $\G$ on some sectors.}
However, even if $\xi$ and $\Gamma$ are real, these complex curves
might not {provide a (real) asymptotic trajectory} if the so called anti-Stokes directions of $\Gamma$ contain the real one.
So even for real analytic vector fields, the theory of multi-summability does not solve the problem we address here.
It should be mentioned that Ecalle proposed in \cite{Eca}
a strategy for a real resummation (see also \cite{Eca-M}). {Our approach circumvent the theory of resurgent functions}.

Still in the real analytic setting, our result has application to tame geometry of trajectories, in the vein of \cite{Rol-San-Sch, LeG-San-Spe1, LeG-San-Spe2, LeG, LeG-Mat-San}. In the last section we prove the following.

\begin{theorem}
\label{cor:ana}
Let $\xi$ be a germ of analytic vector field at
$a\in \R^{m}$, and $\G$ a formal curve at $a$, invariant by $\xi$, 
that is not contained in the singular locus of $\xi$. 
Then, there exists a germ of trajectory $\g$ of $\xi$ asymptotic to $\Gamma$ and subanalytically non-oscillating.
\end{theorem}

We now outline the main steps of the proofs of Theorem~\ref{thm:1} and Theorem~\ref{thm:2}.

After a sequence of admissible transformations (blow-ups and ramifications), 
we reduce the vector field $\xi$ in a neighborhood of the formal curve $\Gamma$, in what we call a Turrittin-Ramis-Sibuya form (TRS for short). 
This reduction is inspired by Turrittin's process \cite{Tur} (see also
Wasov \cite{Was} or \cite{Bal}), developped for linear meromorphic differential equations over $\mathbb C$. 
Ramis and Sibuya used such normal forms for their analysis of multi-summability of formal solutions of (non linear) meromorphic ODE's in \cite{Ram-S}.
L\'opez-Hernanz, Rib\'on, Sanz S\'anchez, Vivas  present in \cite{LRSV} a reduction of the same nature for 
germs of holomorphic diffeomorphisms.
For our purpose, it is necessary to retain the real structure, so we build over Barkatou, Carnicero, Sanz S\'anchez \cite{Bar-Car-San} which gives a real reduction for linear formal meromorphic differential systems.
Our reduction to real (TRS)-form {for a $C^\infty$ vector field along a formal invariant curve} is presented in section~\ref{sec:TRS}. 

Once the vector field is reduced, our general strategy to construct the curve $\gamma$ is to work inductively on the dimension of the ambient space, 
by restriction to a center manifold, until this dimension drops to $1$, {or $\G$ is tangent to a non-zero eigenvalue} (see Lemma \ref{prop:key} in section~\ref{sec:key}).
{In ambient dimension $m=2$, a vector field in TRS-form is either hyperbolic or has a center manifold of dimension 1, so the result is consequence of the classical theory of invariant manifolds (for instance in \cite{Hir-Pug-Shu}).} 
In higher dimension, this approach leads to two main difficulties. At first, no smooth center manifold 
exists in general (see \cite{vStr}), so we have to consider vector fields of finite differentiability class. 

A more serious obstruction appears when the center manifold 
is the full ambient space (all eigenvalues have zero real part) since {in this case the induction is interrupted}.
{Thus, we need new arguments to treat the so called dominant rotation case, excluded in Lemma \ref{prop:key}. That is, when $\xi$ has only eigenvalues with pure imaginary 
initial part.}
In dimension $3$, this situation corresponds to
the rotation case of \cite[IV (2.2) p. 134]{Bon}, treated appart by Bonckaert and Dumortier in \cite{Bon-D}.
{The strategy proposed in that paper} consists of building an invariant slow manifold -- i.e., tangent to the kernel of the linear part of $\xi$ --,
in a similar way center manifolds are constructed in the general theory.
In higher dimension, different rotations with different orders might compete with many real slow directions
of different orders also, and the calculation of the needed estimates seems impracticable.

To deal with it, we introduce in section~\ref{sec:straighten} a special kind of transformations we call straighteners.
{They act} as a direct sum of plane rotations over a fibration transverse to $\G$
so to annihilate the spiraling effect
induced by the pure imaginary eigenvalues.
These transformations are strongly irregular and do not admit even a continuous extension at the singular point.
 {However, in a neighborhood of each half branch of $\G$, 
 $\xi$ has a lift of any finite differentiability class by the convenient straightener, 
 and this lift has no more dominant rotation,
up to first reduce to a stronger TRS-form (the class depends on the strength).
{From here, the induction can be continued. This way, we produce trajectories with high but finite contact order with $\G$.} 

Our final argument to get trajectories asymptotic to $\G$ is based on the existence of the so called accompanying curves in the center manifold, that permits to show that all trajectories with a sufficiently high contact with $\Gamma$ have flat contact the ones with the others (cf. points (ii) in Lemma \ref{prop:key} and Lemma \ref{prop:unlace}). We recall in section~\ref{sec:key} the basics about accompanying curves, deduced from a fine treatment of the principle of reduction to a center manifold by Carr \cite{Carr}. This approach has already been used by Cano, Moussu and Sanz in \cite{Can-Mou-San2} for three-dimensional analytic vector fields. {In this paper, we adapt Carr's construction in order to obtain, moreover, the manifold $S^+$ of Theorem~\ref{thm:2}. A precise statement and the} details are discussed in section~\ref{sec:final}.

 It would be interesting to extend our result to more general type of series. {For instance,} allowing real exponents, or for more general formal transseries.} 
The later needs an extended notion of been asymptotic, a {question that has been considered by} vdHoeven in \cite{vdH}
and by Ashenbrener, vdDries, vdHoeven in \cite{Ash-vdD-vdH} in the context of polynomial ODE's over Hardy fields. 

\subsection{Notations}
Consider a $C^{\infty}$ manifold $M$ of dimension $m$ and $a\in M$. 
We will often put $m=1+n$ and have a local system of coordinates $(x,y_1,\dots,y_n):M\to\mathbb R^{1+n}$, centered   at $a$,
with a distinguished first coordinate. For short, we use a bold letter to refer to the tuple whose 
components are written with the same letter and subscripts: $\yy=(y_1,\dots,y_n)$. We also use subscripts to indicate coordinates or tuples of coordinates of a given object, e.g.,
a parametrized curve $\gamma:\mathbb R \to M$ might be written with no other precision as $(\gamma_x,\gamma_{\yy})$. 

The differential of a map $f$ at a point $a$ is 
written $df(a)\in (T_aM)^*$, where $(T_aM)^*$ is the cotangent space of $M$ at $a$; $a$ might be omitted depending on the context. 
We use symbolic powers for diagonal $k$-tuples, that is, $d^kf(a)(v^{(k)})$ is the value of the $k$-th differential of $f$ over the $k$-tuple of direction $(v,v,\dots,v)\in (T_aM)^k$. 
Given coordinates $(x,\yy)$, the dual basis of $(dx,d\yy)$ is denoted $(\partial_x,\partial_{\yy})=(\partial_x,\partial_{y_1},\dots,\partial_{y_n})$.
We write the action of derivations as a product or with parenthesis, so $\partial_{\yy}(f)=df(\partial_{\yy})=(\partial_{y_1}f,\dots,\partial_{y_n}f)$, not to be mixed up with the composition, e.g., $\xi\circ\gamma$
is the value of $\xi\in\text{Der}_{\mathbb R}(C^{\infty}(M))$ over the parametrized curve $\gamma$; for the later we might also use restriction notations, e.g., $\xi_{|a}$ is the value of $\xi$ at $a$.  
We use a generic $\cdot$ symbol to indicate a dot product on diverse tuples (matrix, vectors, $\dots$). Together with the bold notations and automatic  definitions by subscript, we get compact expressions like $\xi = \xi_x \partial_{x}+\xi_{\yy}\cdot \partial_{\yy}$, where $\xi_x = \xi(x)$ and $\xi_{\yy}=\xi(\yy)$ are implicitely defined once $\xi\in\text{Der}_{\mathbb R}(C^{\infty}(M))$ and $(x,\yy)$ are given. 

We use multi-index for higher order derivatives,
this means, if $\alpha=(\alpha_0,\dots,\alpha_n)\in\mathbb N^{1+n}$ we set
$|\alpha|=\sum_{j=0}^{n} \alpha_j$ and $(x,\yy)^{\alpha}=x^{\alpha_0}y_1^{\alpha_1}\dots y_n^{\alpha_n}$ and then
$$\partial^{|\alpha|}_{(x,\yy)^{\alpha}} = (\partial_x)^{\alpha_0}(\partial_{y_1})^{\alpha_1}\dots(\partial_{y_n})^{\alpha_n}.$$ The jet $j_k f$ at $(0,\bm 0)$ of order $k$ of a function $f$ is the polynomial $$j_kf(0,\bm 0)(x,\yy) = \sum_{j\le k}\frac{1}{j!} d^{j}f(0,\bm 0)((x,\yy)^{(j)}) = \sum_{|\alpha|\le k} \frac{1}{|\alpha|!}(\partial^{|\alpha|}_{(x,\yy)^{\alpha}}{f})(0,\bm 0) (x,\yy)^{\alpha},$$
and the Taylor expansion of $f$ is written $\wh f$. 
We identify polynomials and polynomial functions, so $\mathbb R[x,\yy]$ is seen as a subset of both formal series $\mathbb R[[x,\yy]]$ and smooth functions $C^{\infty}(\mathbb R^{1+n})$. We write $\mathbb R_k[x]$ for the set of
polynomials of (total) degree at most $k$.

We use Landau notations $o$ and $O$ in the
$C^{\infty}$ context, locally at $a$: $f=O(g)$ (resp. $f=o(g)$) if 
there exists a bounded function $h$ (resp. $h$ tends to $0$) such that $f=gh$ in a neighborhood of $a$. 
Notice that when $g$ is a power of a coordinate function, say $g=x^k$, $f=O(x^k)$ (resp $f=o(x^k)$) implies 
$f$ is divisible by $x^k$, that is, $f=x^kh$ with a $C^{\infty}$ function $h$. We use Landau
notations also to compare formal series with powers of coordinates,
which of course implies divisibility in the ring $\mathbb R[[x,\yy]]$.

Given a ring $R$, $\MM_n(R)$ and $\GL_n(R)$ refer respectively to $n\times n$ matrices with coefficients in $R$ and invertible such matrices. We write $I_n$ for the identity matrix in $\GL_n(R)$. Let us introduce some useful notations concerning real and complex matrices. 
Recall $\Theta: \mathbb C\ni a+ib \mapsto a I_2+ b J_2\in \MM_2(\mathbb R)$ is an isomorphism between $\mathbb C$ and the subspace of $\MM_2(\mathbb R)$ spanned
by $I_2, J_2$, with $$I_2=\left(\begin{array}{cc}1 & 0 \\ 0 & 1\end{array}\right) \text{ and } J_2 = \left(\begin{array}{cc}0 & -1 \\ 1 & 0\end{array}\right).$$
We extend $\Theta$ first to formal series by setting, if $h(x)\in\mathbb C[[x]]$: $\Theta(h(x))=\text{Re}(h(x))I_2+\text{Im}(h(x))J_2$.
Then, we let $\Theta$ act on each coefficient of a given matrix $M\in\MM_m(\mathbb C[[x]])$ to define $\Theta(M)$, a matrix \emph{a priori} in $\MM_m(\MM_2(\mathbb R[[x]]))$,
space we identify with $\MM_{2m}(\mathbb R[[x]])$. This way, for each $m\ge 1$, $\Theta$ defines an injective morphism
of $\mathbb R$-algebras between $\MM_m(\mathbb C[[x]])$ and $\MM_{2m}(\mathbb R[[x]])$. 

Since we work with block shaped matrices, it will be convenient to denote $M\oplus N$ the matrix 
$$
M\oplus N=\left(\begin{array}{cc}
M & 0 \\
0 & N
\end{array}\right) \in \mathcal M_{m+n}(R),
$$
whenever $M\in \mathcal M_m(R)$ and $N\in \mathcal M_n(R)$.
If $D$ is diagonal by blocks, of the form $D = \Theta(c_1 I_{n_{1}}\oplus\dots\oplus c_{k} I_{n_{k}}) \oplus d_1 I_{m_1}\oplus \dots\oplus d_{k'} I_{m_{k'}}$,  we say that $C$ has a block structure compatible with $D$ (or that $C$ is compatible with $D$ for short) whenever $C = \Theta(C_1\oplus C_2\oplus\dots\oplus C_k)\oplus E_1\oplus\dots\oplus E_{k'}$
where for all $j$, $C_j\in\mathcal M_{n_j}(R)$ and $E_j\in\mathcal M_{m_j}(R)$. If $C$ is compatible with $D$, then $[D,C]=0$.

\section{Reduction to Turrittin-Ramis-Sibuya form}\label{sec:TRS}

In this section we give a procedure to transform a $C^{\infty}$ vector field along an invariant formal curve to another one with a useful expression in local coordinates.
In the first subsection, we summarize the results of \cite{Bar-Car-San} we base our reduction on.
In a second subsection, we introduce the transformations that are admissible for a couple formed by a vector field and a non-singular invariant curve. 
At last, we  explain how to reduce a given vector field in a ``neighborhood'' of an invariant formal curve (Theorem \ref{th:TRS-vf}).

\subsection{Real Turrittin's Theorem for linear systems of ODEs}
The reduction we are looking for is based mainly in a result by Barkatou, Carnicero and Sanz \cite{Bar-Car-San}, what we discuss briefly in this paragraph. It consists of a version of a classical Turrittin's result on normal forms of formal meromorphic linear systems of ODEs (see \cite{Tur, Was, Bal}) when the base field of coefficients is $\R$. 

Consider a formal linear system of $n$ ODEs of the form
$$
x^{p+1}\bm{y}'=A(x)\cdot\bm{y},
$$ 
where $\bm y =(y_1,...,y_n)\in\mathbb R^{n}$, the apostrophe denotes the derivation with respect $x$, $p$ is an integer and $A(x)\in \MM_{n}(\R[[x]])$, with $A(0)\neq 0$.
The system is {\em singular} when $p\ge 0$ and in this case the number $p$ is called the {\em Poincaré rank} of the system. 

The reduction is obtained by applying to the system certain transformations of the following kind.
\begin{enumerate}
\item {\em Gauge transformations.-}  
If $T(x)\in \GL_n(\R[[x]][x^{-1}])$, the change of variables $\bm y=T(x)\cdot \bm z$ gives rise to a bijection $\Psi_{T(x)}$ between the whole family of systems, called a {\em gauge transformation}. Explicitly, it maps the system $x^{p+1}\bm y'=A(x)\cdot \bm y$ to the system $x^{\wt{p}+1}\bm z'=B(x)\cdot\bm z$ where
$$
x^{-(\wt{p}+1)}B(x)=x^{-(p+1)}T(x)^{-1}A(x)T(x)-T(x)^{-1}T'(x).
$$
We shall consider the following two types of gauge transformations:
\begin{enumerate}
\item {\em Regular polynomial.-} A transformation $\Psi_{P(x)}$ where $P(x)\in\MM_n(\R[x])$ is a polynomial matrix and $P(0)\in\GL_n(\mathbb R)$.
\item {\em Diagonal Monomial.-} A transformation $\Psi_{T(x)}$ where $T(x)$ is diagonal of the form $T(x)=\text{diag}\,(x^{k_1},x^{k_2},...,x^{k_n})$ for some non-negative integers $k_1,...,k_n$, not all of them equal to zero. 
\end{enumerate}
\item {\em Ramification of order $r\in\N_{>1}$.-} Denoted by $R_r$, it corresponds to the change of the independent variable $x=z^r$. It transforms a system $x^{p+1}\bm y'=A(x)\cdot \bm y$ into the system (re-written with the same variable $x$) 
$$
x^{pr+1}\bm y'=r^{-1}A(x^r)\cdot \bm y.
$$
\end{enumerate}
Given a system (S), a transformation 
is called \emph{admissible} for (S) if it is either a gauge transformation of type (a) or (b) above and $T^{-1}AT-x^{p+1}T^{-1}T'$ belongs to $\MM_n(\mathbb R[[x]])$ (this is always the case if $\Psi_{T(x)}$ is regular polynomial) or a ramification $R_r$ and (S) is singular ($p\ge 0$). By extension, a composition of such transformations is \emph{admissible} for a system $(S)$ if 
each transformation is admissible for the system it is applied to.

To introduce the main result of \cite{Bar-Car-San}, we need the following.
\begin{definition}\label{def:TRS-form-lineal}~
Let $q$ be a non-negative integer. A singular system is said to be in {\em Turrittin-Ramis-Sibuya form of Poincaré rank $q$} (or (TRS)$^q$-form) if it is written as
\begin{equation}\tag{(TRS)$^q$-form}		x^{q+1}\bm y'=(D(x)+x^qC+x^{q+1}V(x))\cdot \bm y,
\end{equation}
where: 
\begin{enumerate}
\item $D\in\MM_n(\mathbb R_{q-1}[x])$ is a polynomial matrix of degree at most $q-1$, and $D(x)=\Theta(D_1(x))\oplus D_2(x)$, with 
$$\begin{array}{c}
D_1(x) = {\rm diag}(c_1(x),\dots,c_{n_1}(x))),\; \forall j=1,\dots,n_1,\; c_j(x)\in\mathbb C_{q-1}[x],\\
D_2(x) = {\rm diag}(d_1(x),\dots,d_{n_2}(x)),\; \forall j=1,\dots,n_2,\; d_j(x)\in \mathbb R_{q-1}[x]\\
\end{array}$$
\item $C\in\mathcal M_n(\mathbb R)$ is compatible with $D(x)$;
\item $\left(D(x)+x^qC\right)_{|x=0}\neq 0$;
\item $V(x)\in \MM_{n}(\mathbb R[[x]])$.
\end{enumerate} 
The matrix $D(x)+x^qC$ is called the {\em principal part} of the system, $D(x)$ and $C$ are called respectively the {\em exponential part} and the {\em residual part} of the system, and $V(x)$ is called the {\em vestigial part} of the system.
\end{definition}

\begin{remark}
Definition~\ref{def:TRS-form-lineal} describes a system of $n=2n_1+n_2$ equations. The splitting between real and complex blocs
is not necessarily unique, but
we will always assume that a given (TRS)$^q$-form of a system has a minimal $n_1$. That is, 
for each $j=1,\dots,n_1$, at least one coefficient of the polynomial $c_j(x)$ is non real. 
\end{remark}

A constant matrix $C\in\MM_{n}(\R)$ will be said to have \emph{good spectrum}
if it has no two eigenvalues (in $\mathbb C$) that differ by a non-zero integer number. 
The main result of \cite{Bar-Car-San} we use goes as follows.

\begin{theorem}[\cite{Bar-Car-San}]\label{th:turrittin-real}\rm
	Consider a singular system $$\hspace{-3cm} (S)\;\;\;\;\;\;\;x^{p+1}\yy'=A(x) \cdot \yy$$ with $A(x)\in\MM_{n}(\R[[x]])$ and $A(0)\ne 0$. Then there exist a ramification $R_r$ and a finite composition $\psi$ of admissible gauge transformations, either regular polynomial or diagonal monomial, such that:
	\begin{enumerate}[(i)]
	\item The composition $\psi\circ R_r$ transforms the system $(S)$ into a system $(\wt{S})$ that is either regular, or in (TRS)$^q$-form for some $q\in\N_{\ge 0}$, with a residual part which has good spectrum.
	\item Let $(\wt S)$ be any system in a (TRS)$^q$-form whose residual part has good spectrum. For any {$N\ge1$}, there exists a regular polynomial gauge transformation $\Psi_{T_N}$, 
	that transforms the system $(\wt{S})$ into another system in (TRS)$^q$-form, with the same principal part as $(\widetilde S)$, and a vestigial part $V$ satisfying $V(x)=O(x^N)$.    %
    \end{enumerate}

\end{theorem}

\subsection{Admissible transformations for vector fields along a formal curve}\label{sec:admissible-for-vf}

Let $a\in X$ be a point in a smooth manifold $X$ of dimension $1+n$, and let 
$(\xi,\G)$ be 
a couple made of a germ $\xi$ of a $C^{\infty}$ vector field at $a$ or a formal vector field $\xi$ at $a$, and a non-singular formal curve $\G$ at $a$, invariant for $\xi$ and not contained in the formal singular locus of $\xi$. We call such couple an \emph{invariant couple}, either \emph{smooth} or \emph{formal} according to the nature of $\xi$, and smooth by default.
We say a system of coordinates $(x,\bm y =(y_1,...,y_n))$ centered at $a$ is {\em adapted to $\G$} if the tangent line of $\G$ is transverse to the hyperplane $x=0$. In such coordinates, $\G$ can be parametrized by $x$. This means that there is a unique $\Gamma_{\yy}=(\Gamma_{y_1},\dots\Gamma_{y_n})\in (x\mathbb{R}[[x]])^n$
such that $\Gamma$ is given by $\bm y=\Gamma_{\yy}(x)$. We also write $\Gamma = (x,\Gamma_{\yy}(x))$.
An adapted system $(x,\yy)$ is said to have \emph{contact order} $m$ with $\Gamma$, for a given $m\in \mathbb N$, if $\ord_x(\Gamma_{\bm y})= m$.

\begin{remark}\label{rk:tangency-order-b}\rm
	If $\xi_{|a}=0$ and $(x,\yy)$ has contact at least $m$ with $\G$, then $$\ord_x(\xi_{\yy}(x,\bm 0))\ge m,$$ where
	 $\xi_{\yy}=(\xi(y_1),\dots,\xi(y_n))$.
	Indeed, $\Gamma$ being invariant, $\Gamma'\wedge\hat\xi\circ\Gamma=0$, which gives, considering the terms in $\partial_x\wedge\partial_{\yy}$:
	$$\hat\xi_{\yy}(\Gamma) - \hat\xi_x(\Gamma)\Gamma_{\yy}'=0.$$ 
	But $\hat \xi_x(\Gamma)=O(x)$ since $\xi_{|a}=0$, and $\Gamma_{\yy}'(x)=O(x^{m-1})$, 
	so $\hat\xi_{\yy}(\Gamma)=O(x^m)$. 
	Now, $\hat\xi_{\yy}(x,\bm 0)= \hat\xi_{\yy}(\Gamma) - \partial_{\yy} \hat\xi_{\yy}(\Gamma)\cdot(\Gamma_{\yy})+o(\Gamma_{\yy})$,
	and since $\Gamma_{\yy}=O(x^m)$, we get $\xi_{\yy}(x,\bm 0) = O(x^m)$ as claimed.
\end{remark}

We define the transformations allowed for an invariant couple $(\xi,\G)$.  
\begin{definition}\label{def:admissible}
	{\em Let $(\xi, \G)$ be a smooth (resp. formal) invariant couple.
		An {\em admissible transformation for $(\xi,\G)$} is a germ of $C^{\infty}$ map $\phi:(Y,b)\to(X,a)$, where $Y$ is a smooth manifold of dimension $1+n$, of one of the following types:
\begin{enumerate}[(i)]
			\item {\em Isomorphism.-} $\phi$ is a germ of $C^{\infty}$ diffeomorphism. 
			\item {\em Blowing-up.-} There exists a germ $(Z,a)\subset(X,a)$ of smooth submanifold, which is (resp. formally) invariant for 
			$\xi$ and not tangent to $\G$ at $a$, 
			such that $\phi$ is the germ at $b$ of the blowing-up $\pi_Z:Y\to X$ with center $Z$ and $b\in\pi_Z^{-1}(a)$ is the point corresponding 
			to the tangent line of $\G$. When $Z=\{a\}$, we say $\pi_Z$ is a punctual blowing-up.
			\item {\em Ramification.-} There exists a system of adapted coordinates $\tau=(x,\yy)$ for $\G$ such that the hyperplane 
			$H=\{x=0\}$ is (resp. formally) invariant for $\xi$, 
			and there exists some $r\in\N_{>0}$ such that $(Y,b)=(\R^{1+n},0)$ and  $\tau\circ\phi=R_r$, where $R_r$ is the map $R_r(x,\yy)=(x^r,\yy)$. 
		\end{enumerate} 
	}
\end{definition}

For each admissible transformation $\phi:(Y,b)\to(X,a)$, the \emph{lift}, or \emph{transformed couple} $\phi^*(\xi,\G)$ of $(\xi,\Gamma)$ by $\phi$
 is the couple $(\wt{\xi},\wt{\G})$, where
$\wt\xi$ is the germ of $C^{\infty}$ (resp. formal) vector field at $b\in Y$ satisfying $\phi_*\wt\xi=\xi$, and $\wt{\G}$ is the non-singular formal curve satisfying $\widehat\phi(\wt\Gamma)=\Gamma$. The invariance conditions ensure that $\wt \xi$ exists as a smooth (resp. formal) vector field, and the condition
on $b$ for the blowing-up ensures that $\wt\G$ is a formal curve at $b$. 
Noticeably, $(\wt \xi, \wt \Gamma)$ is an invariant couple again.
Iterating, the lift of $(\xi,\Gamma)$ by a finite composition of admissible transformations $\psi=\phi_r\circ\cdots\circ\phi_1$, refers to 
$\psi^*(\xi,\G)=\phi_1^*\phi_{2}^*\cdots\phi_r^*(\xi,\G)$. 

The (TRS)-form we provide below for an invariant couple $(\xi,\Gamma)$ is a practical expression 
of $\xi$ in some coordinates adapted to $\Gamma$, so we often need to reason with particular coordinate
systems. For this, we list below the coordinate systems and change of coordinates we use
and the effect of the admissible transformations on the coordinates of $(\xi,\Gamma)$.

\begin{definition}\label{def:admissible:coordinates}\rm

Let $(\xi,\Gamma)$ be a a smooth (resp. formal) invariant couple, and $(x,\bm y)$ be an adapted coordinate system.
An {\em admissible coordinate transformation for $(\xi,\G, (x,\bm y))$} is a germ of $C^{\infty}$ map 
		$\phi:(\mathbb R^{1+n},0)\ni (x, \yy)\mapsto (\wt x,\wt \yy) \in (\mathbb R^{1+n},0)$,
	(resp. a formal map {$(\wt x, \bm{\wt y})=\phi(x,\bm y)\in \mathbb R[[x,{\yy}]]^{1+n}$})
	of the following types. For each type, we give the
		expression of the transformed couple $(\wt\xi,\wt\Gamma)=\phi^*(\xi,\G)$ in the coordinate system $(\wt x, \wt {\bm y})$.
\begin{enumerate}
\item {\em Affine polynomial.}  It regroups two types of transformations:  
	\begin{enumerate}
	\item {\em Polynomial translation.-} A map of the form
	$$
	(x,\yy)=T_{\bm\beta}(\wt x,\wt \yy):=
	(\wt x,\bm\beta(\wt x)+\wt\yy),
	$$
	where $\bm\beta(\wt x)\in(x\R[\wt x])^{n}$. In coordinates, we get:
	$$\begin{array}{c}
	\wt \xi = \xi_x\circ\phi\; \partial_{\wt x}+(\xi_{\yy}\circ\phi\;-(\xi_x\circ\phi)\; \bm\beta'(\wt x))\cdot \partial_{\wt{\bm y}}\\
	\wt \Gamma = (\wt x, \Gamma_{\bm y}(\wt x) -\bm\beta(\wt x))
	\end{array}$$ 
	\item {\em Polynomial regular.-} Any map of the form
	$$
	(x,\yy)=\Psi_{P}(\wt x,\wt \yy):=(\wt x, P(\wt x)\cdot \wt\yy),
	$$
	where $P(\wt x)\in\MM_{n}(\R[\wt x])$ and $P(0)\in \GL_n(\mathbb R)$. In coordinates, we get:
	$$\begin{array}{c}
	\wt \xi = \xi_x\circ\phi\; \partial_{\wt x}+(P^{-1}(\wt x) \cdot \xi_{\yy}\circ\phi\;-(\xi_x\circ\phi)\,{P^{-1}(\wt x) \cdot} P'(\wt x)\cdot\wt\yy)\cdot \partial_{\wt{\bm y}}\\
	\wt \Gamma = (\wt x, P^{-1}(\wt x) \cdot \Gamma_{\bm y}(\wt x))
	\end{array}$$ 
	($P(0)\in \GL_n(\mathbb R)$ implies $P^{-1}$ exists in $\MM_n(C^{\infty}(\mathbb R))$ and in $\MM_n(\mathbb R[[x]])$).
	\end{enumerate}
\item {\em Diagonal monomial.-} A map $\phi$ of the form $(x,\bm y) = (\wt x, ((\wt x I_k) \oplus I_{n-k}) \cdot \wt \yy)$, with $1\le k \le n$, admissible if
$(x,\yy)$ has contact order at least $2$ with $\Gamma$ and the center $\{x=0,\yy_1=\bm 0\}$ is invariant by $\xi$, where $\yy_1=(y_1,\dots,y_k)$.
In this case, writing also $\yy_2=(y_{k+1},\dots,y_{n}),\; \widetilde\yy_1=(\wt y_1,\dots,\wt y_k),\; \wt\yy_2=(\wt y_{k+1},\dots,\wt y_n)$,
$$\begin{array}{c}
\wt \xi=(\xi_x\circ\phi)\; \partial_{\wt x}+\frac{1}{\wt x}\, \left(\xi_{\yy_1}\circ\phi - (\xi_x\circ\phi)\;\wt\yy_1\right) \cdot \partial_{\wt\yy_1} + (\xi_{\yy_2}\circ\phi)\cdot\partial_{\wt\yy_2}\\
\wt \Gamma = (\widetilde x, \frac{1}{\wt x}\Gamma_{\yy_1}(\wt x),\Gamma_{\yy_2}(\wt x))
\end{array}$$
When $k=n$, we say that the transformation is \emph{full diagonal monomial}. 
\item \emph{Ramifications.-} A map $\phi$ of the form $(x,\yy)=(\wt{x}^r, \wt \yy)$, where $r\in\mathbb N_{>1}$, admissible if the hypersurface $\{x=0\}$ is invariant by $\xi$.
In this case we have
$$\begin{array}{c}
\wt \xi = \frac{1}{r}\; \wt x^{1-r}(\xi_x\circ\phi) \; \partial_{\wt x} + (\xi_{\yy}\circ\phi) \cdot \partial_{\wt \yy} \\
\wt \Gamma = (\wt x, \Gamma_{\yy}(\wt x^r )). 
\end{array}$$
\end{enumerate} 
\end{definition}


The different notions of admissible transformations so introduced are closely linked.
It is clear that an admissible coordinate transformation for $(\xi,\Gamma,(x,\yy))$ is the expression in adapted coordinates 
of an admissible transformation of $(\xi,\Gamma)$. Affine polynomial transformations are isomorphisms.
A diagonal monomial transformation is the expression in coordinates of the blowing-up of the invariant center $\{(x,\yy_1)=0\}$ 
and the contact order condition of $(x,\yy)$ with $\G$ implies the point $(\wt x,\wt\yy) =0$ corresponds to the tangent line of $\Gamma$.
The definition of the ramification as a coordinate transformation coincide with the one as admissible for the couple $(\xi,\Gamma)$.

On the another hand, the gauge transformations and ramifications that are admissible for a formal linear system 
$x^{p+1}\yy'=A(x)\cdot\yy$ correspond to compositions of 
admissible coordinate transformations 
for $(\xi,\Gamma,(x,\yy))$ with 
$\xi = x^{p+1}\partial_x+(A(x)\cdot\yy)\cdot\partial_{\yy}$ and $\Gamma = (x,\bm{0})$. 
This is checked directly, except for diagonal monomial transformations which need the following.
\begin{proposition}\rm
Let $\Psi_T$ be a diagonal monomial transformation admissible for the system $(S):\; x^{p+1}\yy'=A(x)\cdot\yy$, let
$(S'):\; x^{q+1}\yy'=B(x)\cdot\yy$ be the transformed system, let
$\xi = x^{p+1}\partial_x +(A(x)\cdot\yy)\cdot\partial_{\yy}$
and $\Gamma=(x,\bm 0)$.
Then there exist a composition of diagonal monomial coordinate transformations $\Phi$ admissible for 
$(\xi,\Gamma,(x,\yy))$
such that 
$$\Phi^*(\xi,\Gamma,(x,\yy))=(\wt x^{q+1}\partial_{\wt x} +(B(\wt x)\cdot\wt{\yy})\cdot\partial_{\wt{\yy}},(\wt x,\bm 0),(\wt x, \wt y)).$$
\end{proposition}

\begin{proof}
We write $T=\text{Diag}(x^{k_1}I_{n_1},\dots,x^{k_m}I_{n_m})$
with $k_1>k_2>\dots>k_m$ and $n_1+\dots+n_m=n$, 
after gathering the exponents which coincide and eventually permuting
the variables. 
Remark that
$$T = T_1^{k_1-k_2}\cdot T_2^{k_2-k_3}{\cdot} \dots {\cdot} T_m^{k_m},$$ where
$T_m=xI_n$, and for $i<m$, 
$T_i = xI_{n_1+\dots +n_i}\oplus I_{n_{i+1}+\dots+n_m}$.
On the level of gauge transformations, this gives $\Psi_{T}=\Psi_{T_m}^{(m)}\circ\dots\circ\Psi_{T_1}^{(k_1-k_2)}$ where powers are compositions.
Now, from their expressions in coordinates, the $\Psi_{T_i}$ act on a system 
the same way the diagonal monomial coordinate transformation 
$\Phi_{T_i}(\wt x, \wt y) = (\wt x, T_i(\wt x)\cdot\wt y)$ act on the corresponding vector field.
So $\Phi = \Phi_{T_m}^{(m)}\circ\dots\circ\Phi_{T_1}^{(k_1-k_2)}$ solves the problem, up to check that it is admissible.

Let us investigate the admissibility condition on a generic case shaped as the $T_i$, say $G=xI_{\ell}\oplus I_{n-\ell}$.
Decompose $A$ by blocks of sizes $\ell$ and $n-\ell$ as $$A=\left(\begin{array}{cc}A_{11} & A_{12}\\ A_{21} & A_{22}\end{array}\right),\; A_{11}\in\mathcal M_{\ell}(\mathbb R[[x]]),\; A_{22}\in\mathcal M_{n-\ell}(\mathbb R[[x]]).$$
Then 
$$G^{-1}\cdot A\cdot G-x^{p+1}G^{-1}\cdot G' = 
\left(\begin{array}{cc}A_{11} -x^pI_{\ell} & x^{-1}A_{12}\\ xA_{21} & A_{22}\end{array}\right),$$
and we see $\Psi_G$ is admissible for (S) if and only if 
$x$ divides $A_{12}$. 

Since the order of the coefficients above the diagonal cannot increase by iterating such sort of transformations,
we get that $\Psi_T$ is admissible for $(S)$ if and only if $\Psi_{T_1}$ is admissible for $(S)$ and 
$\Psi_{\check{T}}$ is admissible for the system transformed by $\Psi_{T_1}$,
where $\check T = T_1^{k_1-k_2-1}\cdot T_2^{k_2-k_3}{\cdot} \dots {\cdot} T_m^{k_m}$. 
So we only need to prove that the admissibility of 
$\Phi_{T_1}$ follows from the admissibility of  $\Psi_{T_1}$.

From the calculation above, $\Psi_{T_1}$ is admissible for
$(S)$ means that $A_{12}(0)=0$. Write $\yy_1=(y_1,\dots,y_{n_1}), \yy_2=(y_{n_1+1},\dots, \yy_n)$, 
so $$\xi = x^{p+1}\partial_x + (A_{11}\cdot \yy_1+A_{12}\cdot \yy_2)\cdot\partial_{\yy_1}+(A_{21}\cdot\yy_1+ A_{22}\cdot \yy_2)\cdot\partial_{\yy_2}.$$
The restriction of $\xi$ to 
the center $\{x=0, \yy_1=0\}$ is given by
$$\xi\circ(0,\bm 0,\yy_2) = (A_{12}(0)\cdot \yy_2)\cdot\partial_{\yy_1}+( A_{22}(0)\cdot \yy_2)\cdot\partial_{\yy_2},$$
and since $A_{12}(0)$ vanishes if $\Psi_{T_1}$ is admissible, we see that this restriction belongs to $\text{Span}(\partial_{\yy_2})=\text{ker}(dx\wedge dy_1\wedge\dots\wedge dy_{n_1})$.
In other terms, $\{x=0, \yy_1=\bm 0\}$ is invariant by $\xi$. Therefore $\Phi_{T_1}$ is admissible as required.
\end{proof}

Full diagonal monomial transformations (punctual blowing-ups) play an important role in
our reduction. We detail here a more or less classical property of divisibility by {the equation of} the exceptional divisor.

\begin{lemma}\label{lem:division}\rm 
Let $(\xi,\Gamma)$ be a smooth invariant couple, $(x,\yy)$ a system of adapted coordinates, $k\in\mathbb N$
and $\phi$ the full diagonal monomial transformation.
Suppose $\xi|_0=0$ and set $(\zeta,\Delta,(x,\zz))=\phi^*(\xi,\Gamma, (x,\yy))$. 
Then $x^{k}$ divides $\widehat{\zeta}_{\zz}$ (resp. $\widehat{\zeta}_x$) if and only if $x^{k}$ divides $\zeta_{\zz}$ (resp. $\zeta_x$) in $C^{\infty}(\mathbb R^{1+n},0)$.
\end{lemma}

\begin{proof}
The direct implication only is not obvious.
We have
$x\zeta_{\zz}(x,\zz)=\xi_{\yy}(x,x\zz)-\xi_x(x,x\zz)\zz$,
and $\zeta_x(x,\zz) = \xi_x(x,x\zz)$.
The proof for both is analogous, we focus on $\zeta_{\zz}$.
We write the Taylor expansion with integral remainder for $\varphi:t\mapsto \xi_{\yy}(t(x,x\zz))-\xi_x(t(x,x\zz))\zz$ between $t=0$ and $t=1$:
$$\begin{array}{l}x\,\zeta_{\zz}(x,\zz) = \varphi(1) = \displaystyle \sum_{m=0}^k\frac{1}{m!} \left(d^{m}\xi_{\yy}(0,0)((x,x\zz)^{(m)})-d^{m}\xi_x(0,0)((x,x\zz)^{(m)})\zz\right) \dots
\\
\hspace{0.2cm} 
+ \displaystyle  \int_0^1\left( d^{k+1}\xi_{\yy}(t(x,x\zz))((x,x\zz)^{(k+1)})-d^{k+1}\xi_x(t(x,x\zz))((x,x\zz)^{(k+1)})\zz\right) \frac{(1-t)^k \;dt}{k!} \\
 = \displaystyle \sum_{m=0}^k\frac{x^m}{m!} \left(d^{m}\xi_{\yy}(0,0)((1,\zz)^{(m)})-d^{m}\xi_x(0,0)((1,\zz)^{(m)})\zz\right)\dots
 \\
 \displaystyle  \hspace{0.2cm}
 +\; x^{k+1}\int_0^1\left( d^{k+1}\xi_{\yy}(t(x,x\zz))((1,\zz)^{(k+1)})-d^{k+1}\xi_x(t(x,x\zz))((1,\zz)^{(k+1)})\zz\right) \frac{(1-t)^k dt}{k!}
\end{array}$$
In the last expression, the integral is $C^{\infty}$ (in terms of $(x,\zz)$) by Leibniz integral rule, 
so $x^{k+1}$ divides the remainder.
The initial sum is a polynomial, of degree $k$ in $x$. The formal divisibility of the left hand side by $x^{k+1}$ implies this polynomial is identically zero. Thus $\zeta_{\zz}(x,\zz)=x^{k} f(x,\zz)$, where $f$ is the parametric integral. 
\end{proof}

\subsection{Reduction of a vector field to Turrittin-Ramis-Sibuya form along an invariant curve}

As in the precedent paragraph, $(\xi,\Gamma)$ is an invariant couple.

\begin{definition}\label{def:TRS-form-vf}\rm
Let $q\in\mathbb N$ be a non-negative integer. We say that the vector field $\xi$ is in {\em Turrittin-Ramis-Sibuya form} ({(TRS)-form} for short) of type $(q,N,M)$ if there exists a system of coordinates $(x,\yy):(X,a)\to\mathbb R^{1+n}$, such that
\begin{equation}\label{eq:TRS-qNM-vf}
\xi=x^e u(x,\yy)\left[
x^{q+1}\partial_x+\left( \left(D(x)+x^qC \right)\cdot\yy+x^{q+1+N}V(x,x^M\yy)\right)\cdot\partial_{\yy}
\right],
\end{equation}
where $e\in\mathbb N$, $u(x,\yy)$ is a germ of $C^{\infty}$ unit, $D(x)$ and $C$ satisfy the same conditions as in Definition~\ref{def:TRS-form-lineal},
and  $V:(\mathbb R^{1+n},0)\to\mathbb R^n$ is a map germ. 
We say that $(x,\yy)$ is a {\em system of (TRS)-coordinates} for $\xi$. 
The terms $D(x)$, $C$ and $V(x,\yy)$ are  respectively called {\em exponential} part, {\em residual} part, and \emph{vestigial} part of the (TRS)-form. 
If $(\xi,\Gamma)$ is an invariant couple and $(x,\yy)$ is system of (TRS) coordinates for $\xi$ that is adapted to $\Gamma$,
we say that the couple $(\xi,\Gamma)$ is in TRS-form.
\end{definition}

The main result in this section establishes that any invariant couple $(\xi,\G)$ can be reduced to (TRS)-form after finitely many admissible transformations.

\begin{theorem}\label{th:TRS-vf}\rm
Let $(\xi,\G)$ be an invariant couple at $(X,a)$, and assume that {$\xi|_a=0$.} Then:
\begin{enumerate}[(i)]
\item There exists $q\in\mathbb N$ and a finite composition $\psi:(Y,b)\to(X,a)$ of admissible transformations for $(\xi,\G)$ 
such that the transformed couple $(\wt{\xi},\wt{\G})=\psi^*(\xi,\G)$ is in (TRS)-form of type $(q,0,0)$, 
with a $C^{\infty}$ vestigial part 
and a residual part with a good spectrum.
Precisely, 
for all system of coordinates $(x,\yy)$ at $(X,a)$ adapted to $(\xi,\G)$, 
there exists a system of (TRS)-coordinates $(\wt x, \wt\yy)$ at $(Y,b)$ for $\psi^*(\xi,\G)$,
and a finite composition of admissible coordinate transformations $\Psi$ 
such that $\Psi^*(\xi,\G,(x,\yy))=(\psi^*(\xi,\G),(\wt x,\wt\yy))$.

\item Let $(\xi,\Gamma)$ be an invariant couple at $(Y,b)$, with (TRS) coordinates $(x,\yy)$
of type $(q,0,0)$, with a $C^{\infty}$ vestigial part 
and a residual part with a good spectrum.
Given $N,M\ge 0$, there exists a finite composition $\psi_{N,M}:(W,c)\to(Y,b)$ of punctual blowing-ups, 
such that $(\wt \xi,\wt \G) : =(\psi_{N,M})^*(\xi,\G)$ is in (TRS)-form of type $(q,N,M)$, with a $C^{\infty}$ vestigial part.
More precisely, there exist a system of (TRS) coordinates $(\wt x,\wt \yy)$ for $(\wt\xi,\wt\Gamma)$, and a finite composition $\Psi_{N,M}$ of admissible coordinate transformations, 
made of affine polynomial and full diagonal monomial transformations, such that
$$(\Psi_{N,M})^*(\xi,\G,(x, \yy)) = (\wt\xi,\wt\G,(\wt x,\wt \yy)),$$
and $(\wt\xi,\wt\G,(\wt x,\wt\yy))$ has the same exponential part than $(\xi,\G,(x,\yy))$
and a residual part with good
spectrum included in $\mathbb R^*_{-}\times i\mathbb R$.

\end{enumerate} 	
\end{theorem}

\begin{remark}\rm
Independently of its use in the present paper, the theorem above has its own interest. 
In particular, it is a precious tool to study the local dynamic of an analytic
vector field.
For this, let us emphasis that, if $\xi$ is analytic, the proof below works
step by step (with some eventual simplifications due to the fact that formal divisibility implies
analytic divisibility), and shows a stronger form of reduction: only analytic 
admissible coordinate transformations are needed (the centers of the blowing-ups are smooth analytic
manifolds).
Actually, we chose to present here a proof that works directly in the analytic framework 
to use Theorem \ref{th:TRS-vf} for other purposes. The proof in the $C^\infty$ case could be shorten a bit 
starting with a system of coordinates $(x,\yy)$ having flat contact 
with the formal curve $\Gamma$ (i.e. $\Gamma_{\yy}=0$).
\end{remark}

\subsection{Proof of Theorem~\ref{th:TRS-vf}}
We fix initial adapted coordinates $(x,\yy)$ at the point $a\in X$, and write 
$$
\xi=
\xi_x\partial_x+\xi_{\yy}\cdot\partial_{\yy} \text{ and } \Gamma=(x,\Gamma_{\yy}(x)).
$$
We consider several admissible coordinate transformations, starting from $(x,\yy)$ to get to a system of (TRS) coordinates.
To lighten the notations, we do not systematically change the name of the different objects after each transformation,
and we recycle old names when it implies no local confusion. 
We divide the proof in several steps.

\subsubsection{Getting the associated system of ODEs}\label{pr:1}
The invariance condition is expressed in terms of the parametrization of $\G$ as  
$$ (\widehat\xi_x\circ\Gamma) \Gamma_{\yy}'=\widehat\xi_{\yy}\circ\Gamma$$
and since $\Gamma$ is not included in the formal singular locus of $\xi$,
it implies $\widehat{\xi}_x\circ\Gamma\neq 0$. Let $m=\text{ord}_x\;\widehat{\xi}_x\circ\Gamma(x)$.

We apply the polynomial translation $T_{j_{m+1}\Gamma_{\yy}(x)}:(x,\yy)\mapsto (x,j_{m+1}\Gamma_{\yy}(x)+\yy)$. 
This way, the transformed coordinates (still denoted $(x,\yy)$) have contact order at least $m+2$ with $\Gamma$.
Notice that this does not affect the order of $\widehat\xi_x\circ\Gamma$.

We now apply $m$ full diagonal monomial transformations. We need to show that it is admissible. 
Inductively, 
after $k<m$ full diagonal monomial transformations, we  
have $\Gamma_{\yy}=O(x^{m+2-k})${and $\xi_x(0,\bm 0)=0$ (since $(\phi^*\xi)_x = \xi_x(x,x\yy)$, where $\phi$ is any such transformation).} 
Specifying  $(\widehat\xi_x\circ\Gamma) \Gamma_{\yy}'=\widehat\xi_{\yy}\circ\Gamma$ at $x=0$ gives $\xi_{\yy}(0,\bm 0)=0$. 
So the origin is invariant by $\xi$ and $m+2-k\ge 2$ implies $\Gamma$ has contact at least $2$ with $(x,\yy)$,
then the next full diagonal transformation is admissible. 
{Notice also that these transformations} do not affect the order of $\widehat\xi_x\circ\G(x)$.

We rename $\zeta$ the vector field before the $m$ diagonal monomial transformations so to 
keep the notation $\xi$ for the new vector field. In particular, $\xi_x(x,\yy)= \zeta_x(x,x^m\yy)$.
We claim $\xi_x = x^m u(x,\yy)$ where $u$ is a $C^{\infty}$ unit. Indeed, 
writing $\zeta_x(x,\yy)=\zeta(x,\bm 0)+O(\yy)$ we get $\xi_x(x,\yy)= \xi_x(x,0)+x^mO(\yy)$,
so $m=\ord_x(\widehat\xi_x\circ\Gamma)=\ord_x(\xi_x(x,\bm 0)+o(x^{m}))$. Then, formally, $\widehat \xi_x(x,\bm 0)=x^m v(x)$ for some $v(x)\in\mathbb R[[x]]$ with $v(0)\neq 0$,
and $x^{m}$ divides $\widehat\xi_x(x,\yy)-\widehat\xi_x(x,\bm 0)$ with a non unit ratio. According to Lemma \ref{lem:division}, this divisibility holds
in the $C^{\infty}$ class, and we get $\xi_x(x,\yy)=x^m u(x,\yy)$, where $u$ is a $C^{\infty}$ unit.

Still from Lemma~\ref{lem:division}, we can now factor out from $\xi_{\yy}$ the larger power $x^e$ with $e\le m$ which divides $\widehat\xi_{\yy}$,
and finally, $\xi$ can be written 
\begin{equation}\label{eq:bar-xi}
\xi = x^e u(x,\yy) \eta,\; \text{ with } \eta = x^{p+1}\partial_x + \eta_{\yy}(x,\yy)\cdot \partial_{\yy},
\end{equation}
where $p\ge -1$, $e+p+1=m$, $\eta_{\yy}$ is $C^{\infty}$ and $x$ does not divide $\widehat\eta$.

To the vector field $\eta$ we associate the system of ordinary differential equations 
$$
x^{p+1}\yy'=\eta_{\yy}(x,\yy).
$$ 
To study its behavior ``along $\G$'', we consider the system of linear ODEs asociated to $(\eta,\G)$, given by
\begin{equation}\label{eq:formal-A}
	x^{p+1}\yy'=\wh{A}(x)\cdot\yy,\;\mbox{ where }\wh{A}(x)=\partial_{\yy}\widehat{\eta}_{\yy}(x,\Gamma_{\yy}(x)).
\end{equation} 

\subsubsection{Proof of item (i)}
The case where $p=-1$ in (\ref{eq:bar-xi}), corresponds to the case where $\eta$ is not singular at the origin. It is easy to handle: 
after a new full diagonal monomial transformation, {which is} admissible for $(\xi,\G)$, we get, renaming the unit,
$$
\xi=x^{e-1}u(x,\yy)\left(x
\partial_x
+(-\yy+O(x))\cdot\partial_{\yy}\right),
$$
which is in (TRS)-form of type $(0,0,0)$.

Let us assume that $p\ge0$, i.e., $\eta$ is singular at the origin.
Notice that $\G$ is invariant for $\wh\eta$ and not contained in its formal singular locus. 
Since the lift of $\xi$ by an admissible coordinate transformation $\varphi$
is given by $\varphi^*\xi = (\varphi_x)^e\; (u\circ\varphi) \; (\varphi^*\eta)$, 
and $\varphi_x$ is a power of $x$, 
it is sufficient to prove Theorem~\ref{th:TRS-vf} for $(\eta,\G)$.

We apply Real Turrittin's Theorem~\ref{th:turrittin-real} (point (i)) to the system (\ref{eq:formal-A}): there exists $r\in\N_{\ge 1}$ and an admissible finite composition $\psi$ of either polynomial regular or diagonal monomial transformations such that $\phi=\psi\circ R_r$ transforms the system (\ref{eq:formal-A}) into a system 
that is either regular -- a case already treated, so we assume the alternative --, or
in (TRS)$^q$-form for some $q\ge0$, say
$$
x^{q+1}\yy'=
\left(D(x)+x^qC+x^{q+1}V(x)\right)\cdot\yy,
$$
where $D,C, V$ satisfy the conclusion of the theorem.

Translated in terms of vector fields, $\phi$ is a composition of admissible coordinate transformations for $(\theta,(x,\bm 0),(x,\yy))$
with $\theta =  x^{p+1}\partial_x+(\widehat A(x)\cdot\yy)\cdot\partial_{\yy}$, and $\phi^*(\theta,(x,\bm 0),(x,\yy))$
is regular or in (TRS)-form of type $(q,0,0)$.
We want to apply the coordinate transformation $\phi$ to $\eta$ also in order to get a vector field in the desired (TRS)-form. We need to prepare first $\eta$ by means of some additional admissible transformations in order $\phi$ is admissible. We use the following lemma.

\begin{lemma}\label{lm:finite-determinancy}\rm
Let $\theta=x^{q+1}\partial_x+(\wh A(x)\cdot \yy) \cdot\partial_{\yy}\in\text{Der}_{\mathbb R}(\mathbb R[[x,\yy]])$, $q \ge 0$, and
$\tau:\mathbb R^{1+n}\ni (\wt x, \wt{\yy}) \mapsto(x,\yy)\in\mathbb R^{1+n}$ be a composition of admissible coordinate transformations for $(\theta,(x,\bm 0),(x,\yy))$.
For all $s\ge 1$, there exists $h(s)$ such that, for all smooth invariant couple $(\xi,\Gamma)$, 
if $$\xi = j_{h(s)}\theta + o(x^{h(s)})\cdot \partial_{\yy} \text{ and } \Gamma_{\yy}=o(x^{h(s)}),$$ 
then $\tau$ is admissible
for $(\xi,\Gamma,(x,\yy))$ and $\tau^*\xi = j_s (\tau^*\theta) + o(\wt x^s)\cdot\partial_{\widetilde\yy}$.
\end{lemma}
\begin{proof}
If the statement of the lemma holds for $\theta$ and $\tau_1$ and for $\tau_1^*\theta$ and $\tau_2$ with respective ``shift functions'' $h_1$, $h_2$, it holds 
for $\theta$ and $\tau_2\circ \tau_1$ with shift function $h_1\circ h_2$. So we simply need to prove it when $\tau$ is an admissible coordinate transformation for $\theta$.
From their expressions in coordinates, the lemma holds if $\tau$ is an affine polynomial transformation with $h(s)=s$.
{If $\tau$ is a} diagonal monomial transformation, it holds with $h(s)=s+1$. Indeed, $\xi = j_{h(s)}\theta + o(x^{h(s)})\cdot \partial_{\yy}$
implies that $\xi\circ (0,\yy) = (A(0)\cdot\yy)\cdot\partial_{\yy}=\theta\circ(0,\yy)$, so the center of the blow-up is invariant 
for $\xi$ as soon as it is invariant for $\theta$, being included in $x=0$. Also,  the condition $\Gamma_{\yy}=o(x^{s+1})$ ensures that the coordinates $(x,\yy)$ have sufficient contact with $\Gamma$, so that $\tau$ is admissible for $(\xi,\G)$.
Finally, the lemma is satisfied if $\tau$ is a ramification with $h(s)=s$ (the hypothesis guarantee $x=0$ is invariant by $\xi$).
\end{proof}

In our situation, we consider the bound $m=h(q+1)$ given by Lemma~\ref{lm:finite-determinancy} for the transformation $\phi$ and the vector field $(\theta,(x,\bm 0),(x,\yy))$. We will apply the polynomial translation $T_{\beta}$ where $\beta(x)=j_{2m}(\G_{\yy}(x))$, 
followed by the composition {$\vp$} of $m$ full diagonal monomial transformations: $\vp(x,\yy)=(x,x^{m}\yy)$. A reasoning identical to those in the first paragraphs of \ref{pr:1} 
shows $\vp\circ T_{\beta}$ is admissible for $(\eta,\Gamma,(x,\yy))$.

We rename $(T_{\beta})^*(\eta,\Gamma,(x,\yy))$ as $(\eta,\Gamma,(x,\yy))$. We have $\Gamma_{\yy}=O(x^{2m+1})$, 
so Remark \ref{rk:tangency-order-b} gives $\eta_{\yy}(x,\bm 0)=O(x^{2m+1})$.
Writing $$\eta = x^{q+1}\partial_x + \left(\eta_{\yy}(x,\bm 0)+\partial_{\yy}\eta_{\yy}(x,\bm 0)\cdot\yy +O(||\yy||^2)\right)\cdot \partial_{\yy},$$
we get 
$$\vp^*\eta(x,\yy) = x^{q+1}\partial_x + \left(\frac{1}{x^{m}}\eta_{\yy}(x,\bm 0)+\partial_{\yy} \eta_{\yy}(x,\bm 0)\cdot\yy -mx^{q}\yy +x^{m}O(||\yy||^2)\right)\cdot \partial_{\yy}.$$
Since $(\vp^*\Gamma)_{\yy}=O(x^{m+1})$, $$j_{m}(\partial_{\yy} \eta_{\yy}(x,0)\cdot\yy)=j_{m}(\partial_{\yy} \eta_{\yy}(x,\Gamma_{\yy}(x))\cdot\yy)=j_{m} (\wh A(x)\cdot\yy),$$ and  $\eta_{\yy}(x,\bm 0)=O(x^{2m+1})$ implies $\frac{1}{x^{m}}\eta_{\yy}(x,0)=O(x^{m+1})$. So finally,
$$
\vp^*{\eta}=j_{m}(\theta)
-mx^{q}\yy\cdot{\partial_{\yy}} +o(x^{m})\cdot\partial_{\yy}.
$$

In accordance with Lemma~\ref{lm:finite-determinancy}, $\phi$ is admissible for $\vp^*\eta+mx^q\yy\cdot\partial_{\yy}$, and 
$\phi^*(\vp^*\eta+mx^q\yy\cdot\partial_{\yy}) = j_{q+1}(\theta) +o(x^{q+1})\partial_{\yy}$.
Now, the radial vector field $\rho = mx^q\yy\cdot{\partial_{\yy}}$ is preserved by any admissible coordinate transformation
that is not a polynomial translation and $\rho(0,\yy)=0$.
Since $\phi$ contains no polynomial translation, 
we get that $\phi$ is also admissible for $\vp^*\eta$ and we have
$$(\vp\circ\phi)^*\eta = j_{q+1}(\theta) - \rho +o(x^{q+1})\cdot\partial_{\yy}.$$
Writing $x^{q+1}W(x,\yy)$ for the factor $o(x^{q+1})$ above and making $j_{q+1}(\theta)$ explicit, this gives:
$$
(\vp\circ\phi)^*\eta =
x^{q+1}{\partial_x}+
\left(\left(D(x)+x^q(C-mI_n)\right)\cdot\yy+x^{q+1}W(x,\yy)\right)\cdot{\partial_{\yy}}.
$$
The later is a TRS-form of type $(q,0,0)$, whose residual part, $C-mI_n$, has good spectrum since $C$ has good spectrum, and whose vestigial part $W(x,\yy)$ is $C^{\infty}$. It finishes the proof of Theorem~\ref{th:TRS-vf}, (i).

\subsubsection{Proof of item (ii) of Theorem~\ref{th:TRS-vf}}
Let $N,M\ge 0$ be two natural numbers as in the statement. 
Let $(\xi,\Gamma,(x,\yy))$ be an invariant couple with TRS coordinates of the form $(q,0,0)$ which satisfies the hypothesis of the theorem, point (ii), with a principal part named $D$, and a vestigial part named $W$.

As in the previous paragraph, we consider the formal system of linear ODEs 
$$
(S)\;\;\; x^{q+1}\yy'=\wh{A}(x)\cdot\yy,
$$
 where $\wh{A}(x):=\partial_{\yy}\xi_{\yy}(x,\Gamma_{\yy}(x))$. 
 System $(S)$ satisfies the hypothesis of Theorem~\ref{th:turrittin-real} (ii), that we apply to height $N+M$:  
 there exists an admissible regular polynomial gauge transformation
 $\Psi_{Q},$ $Q(x)\in\MM_{n}(\R[x])$ with $Q(0)=I_n$, that transforms $(S)$ into a system of the form
$$
(\wt{S})\;\;\;\;\;x^{q+1}\yy'=(D(x)+x^qC+x^{q+1+N+M}V(x)))\cdot \yy,
$$
with $D$, $C$, $V$ as in the thesis of the theorem.
We let $m=h(q+1+N+M)$ be the integer given by Lemma~\ref{lm:finite-determinancy}, applied to $\Psi_Q$ and the formal vector field $\theta = x^{q+1}\partial_x+(\wh A(x)\cdot\yy)\cdot \partial_{\yy}$.

We consider two positive integer numbers $\ell,\ell'$ satisfying the following conditions
$$
\begin{array}{l}
\ell\ge\max\{m-q,\,\max(\text{Spec(C)})-M+1\}\\
\ell'\ge\max\{\ell+m,\ell+M+1\}
\end{array}
$$

We apply the polynomial translation $T_\beta$ where $\beta(x)=j_{\ell'}(\G_{\yy}(x))$, so we assume the coordinates $(x,\yy)$ have contact order at least $\ell'+1$ with $\G$. We then apply the composition $\phi$ of $\ell$ full diagonal monomial transformations, this is, $\phi(x,\wt{\yy})=(x,x^\ell\wt{\yy}).$
As before, $\phi$ is an admissible transformation for $(\xi,\G,(x,\yy))$ since $\ell'>\ell$.
We write $(\wt\xi,\wt\G,(x,\wt\yy)):=\phi^*(\xi,\G,(x,\yy))$.
We get the following expression for $\wt{\xi}$:
$$
\wt{\xi}=x^{q+1}\partial_x+
\left(\frac{1}{x^{\ell}}\,\xi_{\yy}(x,\bm 0)+(\partial_{\yy}\xi_{\yy}(x,\bm 0)-\ell x^qI_n)\cdot\wt\yy+\frac{1}{x^{\ell}}\,\bm \chi(x,x^\ell\wt{\yy})\right)\cdot\partial_{\wt{\yy}},
$$
where $\bm \chi(x,\yy)=\xi_{\yy}(x,\yy)-\xi_{\yy}(x,\bm 0)-\partial_{\yy}\xi_{\yy}(x,\bm 0)\cdot\yy$.
The particular form of $\xi_{\yy}$ shows that we also have 
$$\bm \chi(x,\yy)=x^{q+1}\left(W(x,\yy)-W(x,\bm 0)-\partial_{\yy}W(x,\bm 0)\cdot\yy\right).$$
So $\bm \chi(x,\yy)=x^{q+1}O(||\yy||^2)$, and then $\frac{1}{x^{\ell}}\bm\chi(x,x^{\ell}\wt\yy)=x^{q+1+\ell}O(||\wt \yy||^2)$.
Since the coordinates $(x,\yy)$ have contact order at least {$\ell'$+1} with $\G$, we deduce $\xi_{\yy}(x,\bm 0)=O(x^{\ell'+1})$
(by Remark~\ref{rk:tangency-order-b}) and hence $j_{\ell'}(\partial_{\yy}\xi_{\yy}(x,\bm 0))=j_{\ell'}(\partial_{\yy}\xi_{\yy}(x,\Gamma_{\yy}(x)))=j_{\ell'}(\wh{A}(x))$. 
Finally, using that $\ell'+1-\ell>m$ and that {$\ell>m-(q+1)$} we have
$$
\wt\xi + \ell x^q \wt\yy\cdot\partial_{\wt\yy} = j_m\theta +o(x^m)\partial_{\wt\yy} \text{ and }\wt\Gamma_{\wt\yy}=o(x^{m}).
$$
From Lemma~\ref{lm:finite-determinancy}, $\Psi_Q$ 
is admissible for $\wt\xi + \ell x^q\wt\yy\cdot\partial_{\yy}$, then for 
$\wt\xi$ (as before, the radial vector field $\ell x^q\wt\yy\cdot\partial_{\wt\yy}$ is invariant by $\Psi_Q$), and we have, re-writing $(x,\yy)=\Psi_Q^*(x,\wt \yy)$,  
$$\Psi_Q^*\wt\xi = x^{q+1}{\partial_x}+
\left(\left(D(x)+x^q(C-\ell I_n)\right)\cdot\yy+x^{q+1+M+N}\bm{\wt\chi}(x,\yy)\right)\cdot{\partial_{\yy}}
$$
for some $C^{\infty}$ map $\bm{\wt\chi}(x,\yy)$.

To conclude, we apply a composition of $M$ full diagonal monomial transformations to $\Psi_Q^*\wt\xi$. It is admissible due to our choice of $\ell'$, because, being $\Psi_Q$ a regular polynomial transformation, we have
$\text{ord}_{x}(\phi^*\wt \Gamma)_{\yy}=\text{ord}_x(\wt\G_{\wt\yy}) \ge \ell'+1-\ell\ge M+2$, and the expression of the
transformed vector field (now denoted $\xi$ again) in terms of the new coordinates (still denoted $(x,\yy)$) is
$$\xi =  x^{q+1}{\partial_x}+
\left(\left(D(x)+x^q(C-(\ell+M) I_n)\right)\cdot\yy+x^{q+1+N}\bm{\wt\chi}(x,x^M\yy)\right)\cdot{\partial_{\yy}}.
$$
It is a (TRS)-form of type $(q,N,M)$, with a $C^{\infty}$ vestigial part $\bm{\wt\chi}$, 
the same exponential part $D$ than the vector field we started with. Moreover, the
residual part $C-(\ell+M)I_n$ has a good spectrum (since $C$ has good 
spectrum) and $\text{Spec}(C)$ is included in $\mathbb R^*_-\times i\mathbb R$ {by the choice of $\ell$}. \hfill{$\square$}

\section{Restricting to a center manifold}\label{sec:key}

In this section, we consider a vector field $\xi$ in (TRS)-form of type $(q,N,0)$, with a vestigial part in a finite differentiability class $C^k$, and with an exponential part whose eigenvalues have no pure imaginary dominant terms. We say that $\xi$ has no dominant rotation.
Being $(x,\yy)$ some (TRS)-coordinates, we prove that $\xi$ admits a trajectory in each half space $\{x>0\}$ and $\{x<0\}$ accumulating to the origin and having a high contact order (related to $N$) with the $x$-axis. We also establish that any two trajectories living in the same half space and accumulating to the origin have flat contact, and we determine the structure of the pencil of all those trajectories. 
We consider finite differentiability classes because our induction involves the restriction to a center manifold that might not be smooth. {Moreover, Lemma \ref{prop:key} below is used in the next section, after a transformation that turns smooth into $C^k$ vector fields.

The main result of this section relies partly on the fact, certainly folklore among specialists, 
that a trajectory $\gamma$ of a vector field {accumulating to a singular point $a$}
and tangent to a local center manifold $W^c$ admits an ``accompanying trajectory'' $\delta$
in the center manifold
that is ``exponentially closed'' to $\g$. The definition of being exponentially closed depends on the parametrization considered for $\g$ and $\delta$. In Carr's book \cite{Carr}, the difference between those trajectories is estimated by a negative exponential in terms of the natural time of the flow of the vector field. We restate below the precise result we need. Namely, if the vector field is a system of ODEs in a privileged coordinate $x$, as happens for a (TRS)-form, the two trajectories $\g$ and $\delta$ have flat contact with respect to $x$, that serves as a common parameter.

We start summarizing Carr's results from \cite[Ch. 2]{Carr}. 
Consider a vector field $\xi$ of class $C^k$ in a neighborhood of $0\in\R^m$, whose linear part, 
in coordinates $(\xx,\ww)\in\R^c\times\R^u=\R^m$ is written as $d\xi(0)=A\oplus B$, 
where $\text{Spec}(A)\subset i\R$ and $\text{Spec}(B)\subset\R^*_+\oplus i\R$. The classical Center Manifold 
Theorem asserts that, in a neighborhood of $0$, there is a subvariety $W^c$, of class $C^k$, 
tangent to $\{\ww=0\}$ and locally invariant for $\xi$ 
(Carr proposes a construction in his book; for the proof of the differentiability class, the reader may consult Kelley \cite{Kel}). 
Fix such a local center manifold $W^c$ and write it as a graph $W^c=\{\ww=h(\xx)\}$, 
where $h$ is $C^k$ in a neighborhood $U_0$ of $0\in\R^c$ and satisfies $h(0)=0,\; dh(0)=0$. Define the projection 
$$
\pi:U_0\times\R^u\to W^c,\;\;\pi(\xx,\ww)=(\xx,h(\xx)). 
$$

\begin{theorem}[Carr]\label{th:Carr}\rm
With the notations above, there is a neighborhood $U$ of $0\in\R^m$ and some constant $\beta>0$, such that the following holds. If $\g:(-\infty,0]\to U$ is a trajectory of $\xi$ such that $\alpha(\g)=0$, {that is, $\lim_{t\to-\infty}\g(t)=0$.} Then there is a unique trajectory  $\delta:(-\infty,0]\to W^c$ of the restriction $\xi|_{W^c}$, called the {\em accompanying trajectory} of $\g$ in $W^c$, such that 
\begin{equation}\label{eq:exp-accompanying}
\|\g(t)-\delta(t)\|=O(e^{\beta t}),\;\;\mbox{ when }t\to-\infty.
\end{equation}
Moreover, there exists a local homeomorphism $\Psi:(U,0)\to(U',0)$, preserving $\ww$ and satisfying the following. Given a trajectory $\g$ in $U$, if $\delta$ is the trajectory in $W^c$ determined by the initial condition $\delta(0)=\pi(\Psi(\g(0)))$, then 
$\alpha(\g)=0$ if and only if $\alpha(\delta)=0$,
and in that case $\delta$ is the accompanying trajectory of $\g$ in $W^c$.
\end{theorem}

\begin{remark}\label{rk:homeo-Carr}
{\em
The homeomorphism $\Psi$ in the statement above is the inverse of the one described by Carr in \cite[p. 22]{Carr} (called $S$ there).  We find useful to use instead the map 
$$
\wt{\Psi}:U\to W^c\times\R^u,\;\;(\xx,\ww)\mapsto(\pi(\Psi(\xx,\ww)),\ww),
$$
that is also a local homeomorphism by the Invariance Domain Theorem.
}
\end{remark}

\begin{proposition}\label{pro:accompanying}\rm
Let $\xi$ be a vector field of class $\CC^k$ in a neighborhood of the origin of $\R^{1+n}$, with coordinates $(x,\yy)$, such that
$$
\xi=x^{q+1}\partial_x+\xi_{\yy}\cdot\partial_\yy,
$$
where $q\ge 1$. Then $\xi$ admits a local $C^k$ center manifold $W^c$, transverse to $\{x=0\}$ and included in a central unstable manifold $W^{cu}$, which satisfies the following. 
\begin{enumerate}
\item If
$\gamma$ is a trajectory parametrized as $(x,\bm\g(x)), x>0$ with $\alpha(\g)=0$ then there exists a trajectory $\delta$ parametrized as $(x,\bm\delta(x)), x>0$, such that $\forall \ell\in\mathbb N,\;||\bm\delta(x)-\bm\g(x)||=o(x^\ell)$.
\item The homeomorphism $\Psi$, associated by Theorem \ref{th:Carr} with a given system of coordinates $(\bm x,\bm w)$ of $W^{cu}$ with $\bm x_1=x$, preserves $x$, i.e., $x\circ \Psi = x$.
\end{enumerate}
\end{proposition}

\begin{proof}
Being $q\ge 1$, $\partial_x\in \text{ker}(d\xi(0))$, so any local center manifold of $\xi$ is transverse to $\{x=0\}$. Let $\g:(-\infty,s)\ni t\mapsto (x(t),\bm\g(x(t)))\in\mathbb R^{1+n}$ be as in the statement,
where $t$ is the natural time associated with $\xi$. Since $\alpha(\g)=0$, $\gamma$ is included in any local center-unstable manifold (that is, invariant and tangent to the sum of the eigenspaces associated with eigenvalues in $\mathbb R_{\le 0}\oplus i\mathbb R$, see \cite{Kel} again). We fix a $C^k$ center-unstable manifold $W^{cu}$, then a $C^k$ center manifold $W^c\subset W^{cu}$ of the restricted vector field $\xi|_{W^{cu}}$. Note that $W^c$ is a local center manifold 
of $\xi$, so proving the proposition inside $W^{cu}$ suffices. We therefore assume $W^{cu}$ and $\mathbb R^{1+n}$ coincide locally, so Theorem~\ref{th:Carr} applies. 
Shrinking our neighborhood of $0$ if needed, we get the accompanying trajectory $\delta:(-\infty,s')\ni t \mapsto(\delta_x(t),\delta_{\yy}(t))\in \mathbb R^{1+n}$ of $\g$ in $W^c$, parametrized by the natural time  (that is, satisfying (\ref{eq:exp-accompanying})).
From $\xi_x=x^{q+1}$ we get that the two functions $x(t)$, $\delta_x(t)$ satisfy the same differential equation 
\begin{equation}\tag{E} \frac{dx}{dt}=x^{q+1}.\end{equation}
{By integration,} the difference of two different solutions of $(E)$ is always larger than a {certain} power of $|t|$ when $t$ goes to $-\infty$. 
But $$|\delta_x(t)-x(t)|\le ||\delta(t)-\gamma(t)||=O(\exp(\beta t)) \text{ as }t\to-\infty,$$ so $\delta_x(t)$ and $x(t)$ cannot be different solutions: $\delta_x(t) = x(t)$. 
Writting now $t(x)$ for the inverse function of $x(t)$, and $\bm\delta(x) = \delta_{\yy}(t(x))$, we see $\delta$ is parametrized
by $(x,\bm\delta(x)), x>0$, and we have
$$||\bm\delta(x)-\bm\g(x)||=||\delta(t(x))-\gamma(t(x))||= O(\exp(\beta t(x)).$$
From $(E)$ we also get $t(x)\sim -\frac{ x^{-q}}{q}$ as ${x\to 0}$, which, replaced in the previous equation, implies point (1).

Following the notations introduced with Theorem \ref{th:Carr}, $\pi\circ\Psi$ maps a point $(x,\bm \gamma(x))$
to $\pi((\delta_x(t(x)), \delta_{\bm y}(t(x)))$. We saw $\delta_x(t(x))=x$, which is point (2).
\end{proof}

The principal result of this section is the forthcoming lemma, which involves the following definitions.
\begin{definition}\label{def:no-dominant-rotation}\rm
Let $q\in\mathbb N_{>1}$, and
$D(x)=\Theta(D_1(x))\oplus D_2(x)\in\MM_n(\mathbb R_{q-1}[x])$, with 
$$\begin{array}{c}
D_1(x) = {\rm diag}(c_1(x),\dots,c_{n_1}(x))),\; \forall j=1,\dots,n_1,\; c_j(x)\in\mathbb C_{q-1}[x],\\
D_2(x) = {\rm diag}(d_1(x),\dots,d_{n_2}(x)),\; \forall j=1,\dots,n_2,\; d_j(x)\in \mathbb R_{q-1}[x].\\
\end{array}$$
We say that $D$ has \emph{no dominant rotation} when $$\forall j=1,\dots,n_1, \; \text{ord}_x\,\text{Re}(c_j(x))\le \text{ord}_x\,\text{Im}(c_j(x)).$$
We also define the {\em unstability index} of $D(x)$ as 
$$
u(D(x)):=2\text{Card}\{j;\;\text{Re}(c_j)>0\}+\text{Card}\{j;\;d_j>0\},
$$
where a real polynomial $P$ satifies $P>0$ iff $P(x)>0$ for $x>0$ sufficiently small.

\end{definition}

\strut

\begin{lemma}\label{prop:key}\rm
Let $\xi$ be a vector field in (TRS)-form of type $(q,N+1,0)$ in the coordinate system $(x,\yy):(X,a)\to\mathbb R^{1+n}$, and suppose:
\begin{enumerate}
\item the exponential part has no dominant rotation;
\item the residual part has spectrum included in $\mathbb R^*_-\oplus i\mathbb R$;
\item the vestigial part is a $C^{n(q+1+N)+1}$ germ.
\end{enumerate} 
Then, 
\begin{enumerate}[(i)]
\item $\xi$ admits a trajectory $(x, \bm\gamma(x)),\; x>0,$ which has contact of order $N+1$ with the $x$-axis: $\bm\gamma(x) = o(x^{N+1})$;
\item If $(x,\bm\zeta(x)), x>0,$ is any trajectory of $\xi$ satisfying $\lim_{x\to 0}\bm\zeta(x)=0$, then $\bm\gamma$ and $\bm\zeta$ have flat contact:
 $\forall k\in\mathbb N,\; \bm\gamma(x)-\bm\zeta(x) = o(x^k)$;
 \item For any neighborhood $U$ of $a$, there is an open neighborhood $V\subset U$ of $a$ and a closed connected topological submanifold $S$ of $V\cap\{x>0\}$ such that:
 \begin{enumerate}
 \item $\dim(S)=1+u(D(x))$;
 \item  $S$ is locally invariant for $\xi$;
 \item For all $b\in V\cap\{x>0\}$, $b\in S$ if and only if the trajectory of $\xi_{|V}$, parametrized as $(x,\bm \gamma(x))$ for $x\in (\alpha,\omega)$, satisfies $\alpha=0$ and 
 $\lim_{x\to 0}\bm\g (x)=0$.
 \end{enumerate}
\end{enumerate} 
\end{lemma}

\begin{proof}
In the definition of a (TRS)-form appears a factor made of a power $x^e$ and a unit $u$.
Since the proposition only concerns the foliation induced by $\xi$ in the half space $x>0$, we might suppose $x^eu(x,\yy)=1$. 
So we suppose $\xi$ is given by:
$$
\xi=
x^{q+1}\partial_x+\left( \left(D(x)+x^qC \right)\cdot\yy+x^{q+N+2}V(x,\yy)\right)\cdot\partial_{\yy}
,$$
$D$, $C$ and $V$ satisfying the hypothesis.

We prove the proposition by induction on the couple $(n,q)$ where $1+n$ is the dimension of the ambient space, and $q$ is the Poincaré rank of $\xi$.
We initialize the induction when $n=0$ or $q=0$ and prove that the case $(n,q)$ follows from a case $(n',q')$ with $n'<n$ and $q'<q$.

{\bf Case $n=0$.}

Here $\xi = x^{q+1}\partial_x$. The curve $x, x>0$ is the only trajectory of $\xi$ included in $x>0$ and the conditions (i)-(iii) are trivial.

{\bf Case $q=0$.}

Here $\xi = x\partial_x + [C\cdot\yy+x^{2+N}V(x, \yy)]\cdot\partial_{\yy}$. 
The lift $\eta$ of $\xi$ by the composition $\yy=x^{N} \zz$ of $N$ full diagonal monomial transformations is given by
$$\eta=  x\partial_x + [(C-N\text{I}_n)\cdot\zz+x^2V(x, x^{N} \zz)]\cdot \partial_{\zz}.$$ 
The origin is an hyperbolic singularity of $\eta$, with 
one positive eigenvalue $1$ associated with $\partial_x$, and $n$ eigenvalues with negative real part (the eigenvalues of $C-N\text{Id}$; recall $\text{Spec}(C)\subset \mathbb R^*_-\oplus i\mathbb R$).
So the unstable manifold of $\eta$ has dimension $1$: it is a trajectory issued from the origin and tangent to $\partial_x$ at $x=0$. In particular, it is not included in $x=0$, and since $\eta(x)\neq 0$ if $x>0$, 
it can be parametrized by $x$. Let $(x,\bm\delta(x))$ be this trajectory. Being tangent to $\partial_x$, we get $\bm\delta(x)=o(x)$.
Then $(x,\bm\gamma(x)) :=(x,x^{N}\bm\delta(x))$ is a trajectory of $\xi$ which satisfies $\bm\gamma(x)=o(x^{N+1})$. This proves (i).

Now, let $\zeta=(x,\bm\zeta(x))$ be any trajectory of $\xi$ with $\lim_{x\to0}\bm\zeta(x)=0$. Such a trajectory is not contained in the stable manifold (which coincides with $x=0$). Thus, since the singularity is hyperbolic and $\zeta$ accumulates to it, $\zeta$ is included in the unstable manifold of $\xi$. This unstable manifold is of dimension one and is made of only one trajectory in $x>0$, which means $\bm\zeta = \bm\gamma$. So conditions (ii) and (iii) are automatically satisfied.

{\bf Induction: $q>0$ and $n>0$.}
We reorder the coordinate functions $\yy$, in such a way that $1,\dots, m$ are the indices such that $\text{ker}(D(0))=\text{Span}(\partial_{\yy_1},\dots,\partial_{\yy_m})$. Write $\zz=(y_1,\dots,y_m)$, $\ww = (y_{m+1},\dots,y_n)$. Recall  ord$_x(D(x)) = 0$ so $m < n$.
Since  $\xi$ has no dominant rotation, the center space (associated with eigenvalues with real part $0$) coincide with the kernel of $d\xi(0)$, that is $\text{Span}(\partial_x, \partial_{\zz})$. So the vectors
$\partial_{\ww}$ are associated with non diagonal elements of $D(0)$.
 Up to permutation, 
$\ww = (\rr,\ss)$, with $\rr = (y_{m+1},\dots,y_{m+1+d})$, $\ss=(y_{m+2+d}, \dots,y_n)$,
and $\partial_{\rr}$ (resp. $\partial_{\ss}$) are associated with positive (resp. negative) diagonal coefficients of $D(0)$.

So $\xi$ admits a $C^{n(q+1+N)+1}$ center unstable manifold $W^{cu}$ and a $C^{n(q+1+N)+1}$ center manifold $W^c\subset W^{cu}$,
that are graphs over $(x,\yy,\rr)$ and $(x,\yy)$ respectively. {The manifold} $W^{cu}$ intervenes only in the proof of point (iii) below, we focus on $W^c$. Call $\bm h$ the $C^{n(q+1+N)+1}$ map that gives $W^c$, that is, 
$$W^c = \{\ww =\bm h(x,\zz)\}.$$

We claim that $\bm h(x,\zz)=x^{q+1+N}\bm g(x,\zz)$ for some $C^{(n-1)(q+1+N)+1}$ map $\bm g$. 
For this, first remark that $W^c \cap \{x=0\}$ {is} a center manifold of the linear vector field $$\xi(0,\yy) = (D(0)\cdot \yy)\cdot\partial_{\yy} = (D_{\ww}(0)\cdot\ww)\cdot\partial_{\ww},$$
where $D_{\ww}$ is defined by $(D(x)\cdot(\zz,\ww))_{\ww} = D_{\ww}(x)\cdot\ww$ (we define $D_\zz$ accordingly, so $D=D_{\zz}\oplus D_{\ww}$). 
This center manifold is clearly unique and given by $\ww=0$, so $W^c \cap \{x=0\}=\{x=0,\; \ww =0\}$. 
This means that $\bm h(0,\zz)=0$, so $\bm h(x,\zz)$ is divisible by $x$ and $\bm h(x,\zz)/x$ is a $C^{n(q+1+N)}$ map. 
Let $$s= \sup \{r\le q+1+N;\; \bm h(x,\zz)/x^r \text{ is }C^1\},$$ and write $\bm h(x,\zz)=x^s \bm g(x,\zz)$.
If $s=q+1+N$, $\bm h$ being $C^{n(q+1+N)+1}$ and divisible by $x^{q+1+N}$, $\bm g$ is $C^{(n-1)(q+1+N)+1}$ and the claim is proven, 
so showing $s=q+1+N$ suffices. 
Suppose $s<q+1+N$, so $\bm g(x,\zz)$ is $C^{n(q+1+N)+1-s}$ and $\bm g(0,\zz)\neq 0$, and let us get a contradiction. 
Since $W^c$ is invariant by $\xi$ and given by $\ww-\bm h(x,\zz)=0$, we have 
${\xi(\ww-\bm h)|_{\ww=\bm h}} = 0$, which gives
$$
\begin{array}{rl}(D_{\ww}+x^q C_{\ww})\bm h = &  \partial_{\zz} \bm h \cdot(D_{\zz}+x^q C_{\zz})\cdot\zz +x^{q+1}\partial_x \bm h \;\dots \\
 & \hspace{0.5cm} +\;x^{q+1+N} \left(\partial_{\zz} h \cdot V_{\zz}(x,\zz,\bm h)-V_{\ww}(x,\zz,\bm h)\right),
 \end{array}
$$
where $C=C_{\zz}\oplus C_{\ww}$ and $V(x,\yy)=(V_{\zz}(x,\yy),V_{\ww}(x,\yy))\in\mathbb R^{m}\times\mathbb R^{n-m}$.
Replacing $\bm h$ by $x^s\bm g$, dividing by $x^s$ and setting $x=0$ in this equation leads to 
$$D_{\ww}(0)\cdot \bm g(0,\zz) = 0,$$ which is a contradiction since $D_{\ww}(0)$ is invertible and $\bm g(0,\zz)\neq 0$.

Now, let $\eta$ be the pullback of $\xi$ by the inclusion $W^c\hookrightarrow \mathbb R^{n+1}$. 
In the coordinate system $(x,\zz)$ of $W^c$, $\eta$ is given by
$$ x^{q+1}\partial_x + \left((D_{\zz}+x^q C_{\zz})\cdot\zz +(D_{\ww}+x^qC_{\ww})\cdot \bm h(x,\zz) + x^{q+1+N}V_{\zz}(x, \zz,\bm h(x,\zz))\right)\cdot\partial_{\zz}.$$
But $$(D_{\ww}(x)+x^qC_{\ww})\cdot \bm h(x,\zz)=  x^{q+1+N}(D_{\ww}(x)+x^qC_{\ww})\cdot \bm g(x,\zz)$$ 
so 
$$\eta = x^{q+1}\partial_x + \left((D_{\zz}(x)+x^q C_{\zz})\cdot\zz + x^{q+1+N}\bm \nu(x,\zz)\right)\cdot \partial_{\zz}$$
for some $C^{(n-1)(q+1+N)+1}$ germ $\bm \nu$. 
Let $v=\text{ord}_x(D_{\zz}(x)+x^qC_{\zz})$, so $1\le v\le q$, since $\text{ord}_x(D_{\zz})\ge 1$ and $C_{\zz}\neq 0$ (recall $0$ is not an eigenvalue of $C$, then neither of $C_{\zz}$). 
The vector field $\eta$ is divisible by $x^v$ and {$\eta/x^{v}$  is in (TRS)-form of type $(q-v,N+1,0)$ in coordinates  $(x,\zz)$} and satisfies:
\begin{enumerate}
\item the exponential part $D_{\zz}(x)$ has no dominant rotation;
\item the residual part $C_{\zz}$ has spectrum included in $\mathbb R^*_-\oplus i\mathbb R$;
\item the vestigial part $\bm \nu(x,\zz)$  is a $C^{(n-1)(q+1+N)+1}$ germ, and $$(n-1)(q+1+N)+1 \ge m(q-v+1+N)+1.$$
\end{enumerate}
So, since $m<n$ and $q-v<q$, the induction hypothesis applies to $\eta/x^v$. 

From this we deduce the three points of the lemma.
For point (i), 
$\eta/x^v$ has a trajectory $(x,\bm\delta(x)), x>0$ such that
$\bm\delta(x)=o(x^{N+1})$. If $\bm\gamma(x) := (\bm\delta(x),\bm h(x,\bm \delta(x)))$, the curve $(x,\bm\gamma(x)), x>0$ is a trajectory of $\xi$, and since 
$\bm \gamma(x)= (\bm \delta(x), x^{q+1+N}\bm g(x,\bm \delta(x)))$, it also verifies $\bm \gamma(x) = o(x^{N+1})$. 

For point (ii), let $\zeta=(x,\bm \zeta(x)), x>0$ be a trajectory of $\xi$ satisfying $\lim_{x\to0}\bm \zeta(x) = 0$. We are in the conditions of Proposition~\ref{pro:accompanying}, so $\zeta$ admits an accompanying trajectory $(x,\bm c(x))$ of $\xi$ contained in 
$W^c$ (that is, $\bm c(x)$ has flat contact with $\bm \zeta(x)$).
As usual, we write $\bm c = (\bm c_{\zz}, \bm c_{\ww})$.
Then $(x,\bm c_{\zz}(x))$ is a trajectory of $\eta/x^v$ and $\lim_{x\to0}\bm c_{\zz}(x)=0$.
From the induction hypothesis, this implies $\bm c_{\zz}(x)$ and $\bm \delta(x)$ have flat contact.
Since $\bm h$ is differentiable, $\bm h(x, \bm c_{\zz}(x))$ and $\bm h(x, \bm \delta(x))$ have flat contact also, 
so $\bm c(x)=(\bm c_{\zz}(x), \bm h(x, \bm c_{\zz}(x)))$ and $\bm \gamma(x)=(\bm \delta(x),\bm h(x,\bm \delta(x)))$ have flat contact. 
Finally, since $\bm \zeta$ has flat contact with $\bm c$ and $\bm c$ 
has flat contact with $\bm \gamma$, we have that $\bm \zeta$ and $\bm \gamma$ have flat contact. {This proves point (ii).}

Let us show point (iii). 
Recall $\yy=(\zz,\rr,\ss)$ where the center manifold $W^c$ is given by $(\rr,\ss)=\bm h(x,\zz)$, and the center-unstable manifold $W^{uc}$ is a graph over $\ss=\bm 0$.
Let $\wt\Psi$ be the homeomorphism introduced in Remark \ref{rk:homeo-Carr}. According to Proposition \ref{pro:accompanying}, point (2), $\wt\Psi$ maps a point
$(x,\zz,\rr,\ss)\in W^{cu}$ to $((x,\zz',\bm h(x,\zz')),\rr)\in W^c\times\mathbb R^d$, where  the trajectory issued from $(x,\zz',\bm h(x,\zz'))$ is the accompanying trajectory of the one issued from $(x,\zz,\rr,\ss)$. We let $\pi$ be the projection onto the first factor of $W^{cu}\times\mathbb R^d$.

Given a neighborhood $U$ of $a$, we fix an open neighborhood $U_0\subset U$ of $a$ in such a way that
$\wt{\Psi}$ restricts to an homeomorphism $\wt\Psi_0$ from $U_0\cap W^{cu}$ onto a product {$U^c_0\times B$} of two connected open sets of $W^c$ and $\R^d$ respectively. 
Let {$V^c\subset U^c_0$} be the open neighborhood of {$a$ in $W^c$}, and $S^c$ be the closed connected topological submanifold of $V^c\cap\{x>0\}$ given by point (iii) for the vector field $\eta/x^v$. Choose an open neighborhood $V\subset U$ of $a$ such that $V\cap W^{cu}=\wt{\Psi}_0^{-1}({V^c\times B})$, and let $S$ be the subset of $V$ given by
$$
S:=\wt{\Psi}_0^{-1}\left({S^c\times B}\right).
$$

Since $\wt\Psi_0$ is an homeomorphism and preserves $x$, $S$ is a closed connected $C^0$ submanifold of $V\cap\{x>0\}$ due to the corresponding property for $S^c$.
The dimension of $S$ is $\dim(S^c)+d$, and according to the inductive hypothesis, $\dim(S^c)=1+u(x^{-v}D_\zz(x))=1+u(D_\zz(x))$. Since $\rr$ is a $d$-tuple, and corresponds to the positive {diagonal elements of $D(0)$}, we conclude $\dim(S)=1+u(D(x))$ (point (iii,a)).

Let $b\in V\cap \{x>0\}$ and call $\gamma$ the trajectory of $\xi_{|V}$ passing through $b$ and parametrized as $(x,\bm\gamma(x)), x\in(\alpha,\omega)$.
If $b\in S$, then $\pi(\wt{\Psi}(b))\in S^c$, so the trajectory $\delta$ passing through $\pi(\wt{\Psi}(b))$ is included in
$S^c$ and has $\alpha$-limit point $0$. By definition of $\wt\Psi_0$, $\delta$ is the accompanying trajectory of $\gamma$, so we deduce $\gamma$ is included in $S$ and $\alpha=0$.
 
The first property shows that $S$ is locally invariant for $\xi$ (point (iii,b)). The second property means that $S$ is composed by trajectories accumulating to $0$ for negative time (point (iii,c) direct implication). 
Suppose now that $b\notin S$. 
If $\alpha=0$ and $\lim_{x\to 0}\bm\g(x)=0$, then $\pi(\wt{\Psi}(\g))$ is a trajectory in $W^c\cap\{x>0\}$ accumulating to $0$, thus contained in $S^c$, contradicting the fact that $b\not\in S$. It proves point (iii,c) reciprocal. \end{proof}

\section{Straightening}\label{sec:straighten}

In this section we prove a proposition analogous to Lemma \ref{prop:key}, but allowing dominant rotation. 
For this, we introduce a particular kind of transformation that we call a straightener.

\begin{definition}\label{def:unlace}
Let $q\ge 0$. A \emph{rotational matrix} of degree $q$ is a polynomial matrix $R\in\MM_n(\mathbb R_{q}[x])$
of the form $$R(x)=\Theta(\text{Diag}(b_1(x),\dots,b_k(x)))\oplus \bm 0_{n-2k},$$ where $\bm 0_{n-2k}\in\MM_{n-2k}(\mathbb R)$ is the null matrix and {$b_j(x)\in i\mathbb R_{q}[x]\setminus\{0\}$ for all $j=1,\dots,k$.}
The \emph{straightener} $U_R$ associated with $R$ is the mapping $U_R(x,\yy)=(x,\Omega_R(x)\cdot \yy)$, where $\yy=(y_1,\dots,y_n)$ and 
$\Omega_R(x)$ is given by:
$$\begin{array}{rrcl} \Omega_R : & \mathbb R^*_+ & \to & \MM_n(\mathbb R) \\
 & x & \mapsto & \displaystyle \exp {\int_x^{+\infty} }\frac{R(t)}{t^{q+2}}\; dt. 

\end{array}$$
The \emph{axis} of the straightener $U_R$ is the linear subspace $y_1=\dots=y_{2k}=0$.
\end{definition}

\begin{remark}\rm
Keeping the notations above and writing, for $j=1,\dots,k$:
$$b_j = i(b_j^0+b_j^1x+\dots+b_j^qx^q) \text{ and }\alpha_j(x) = \frac{b_i^0}{(q+1)x^{q+1}}+\frac{b_i^1}{qx^{q}}+\dots+\frac{b_i^{q}}{x},$$
$\Omega_R$ is given by 
$$\Omega_R = \left(\cos(\alpha_1)I_2+\sin(\alpha_1)J_2\right)\oplus\dots\oplus\left(\cos(\alpha_k)I_2+\sin(\alpha_k)J_2\right)\oplus I_{n-2k}.$$
So the {map} $U_R$ acts on the fibers of $(x,\yy)\mapsto (x,y_{k+1},\dots,y_{n})$ as a direct sum of plane rotations 
of angles $-\alpha_j(x)$, that are unbounded as $x$ goes to $0$.
In particular, say $n=2$ and $k=1$ for simplicity, it interlaces the ``horizontal" lines (the fibers of $\yy$) the ones with the others.
We chose to call it ``straightener''
however,
because we apply it to curves that are already interlaced, but the other way round;
the straightener mapping, by interlacing regular curves, will unlace our curves of interest.
For instance, the vector field $x^{q+2} \partial_x +(R(x)\cdot\yy)\cdot\partial_{\yy}$ have ``spiraling'' 
trajectories $(x,\yy(x))$ given by $$\yy(x) = \displaystyle\left( \exp {\int_x^{+\infty} } \frac{R(t)}{t^{q+2}}\; dt\right) \cdot \yy_0,\; \yy_0\in\mathbb R^2.$$
The lifts $(x,\zz(x))$ of these trajectories by the straightener $(x,\yy)=U_R(x,\zz)$ are the horizontal lines $\zz(x)={\zz_0}$. In this way, $U_R$
straightens the coiling induced by the rotational part $R$ of the vector field.
\end{remark}

We will use the following properties of straighteners.

\begin{lemma}\label{lem:omega}\rm
Let $n\ge 2$, $q\ge 1$, $M\ge 1$, $R$ be a rotational matrix of degree $q-1$, and let $\Omega_R$, $U_R$ be defined as in Definition \ref{def:unlace}. Then:
\begin{enumerate}
\item $U_R$ is a $C^{\infty}$ diffeomorphism of $\mathbb R_+^*\times\mathbb R^n$ which coincide with the identity in restriction to its axis (and in particular to the $x$-axis);
\item $\text{ord}_x$ is invariant by $U_R$: the two curves $(x,\bm\gamma(x))$ and $(x,\bm\delta(x))$
have contact order {at least} $N$ (i.e., $\|\bm\gamma(x)-\bm\delta(x)\|=O(x^N)$) if and only if 
the two curves $(U_R)^*(x,\bm\gamma(x))$ and $(U_R)^*(x,\bm\delta(x))$ have contact of order $N$.
In particular, the contact of a curve with the $x$-axis is invariant by $U_R$;
\item $\Omega_R$ satisfies the differential equation $x^{q+1}\Omega_R' = R\cdot \Omega_R$;
\item the map $x \mapsto x^M \Omega_R(x)$ admits a $C^{\lfloor\frac{M}{q+1}\rfloor}$ extension on $\mathbb R_+$; 
 \item if $C$ is compatible with $R$, then $C$ commutes with $\Omega_R$ and $\Omega_R^{-1}$.
\end{enumerate}
 \end{lemma}

\begin{proof}
\begin{enumerate}
\item $U_R$ admits $U_{-R}$ as a reciprocal, and is clearly smooth; the expression of $U_R$ in restriction to its axis is the identity.
\item For any fixed $x>0$, the map $\yy\mapsto U_R(x,\yy)$ is an isometry (w.r.t. euclidean distance).
\item From the shape of $R$, we notice that $\int_x^{+\infty} \frac{R(t)}{t^{q+1}}\; dt$ commutes with its derivative. The result follows classically.
\item By induction on $\lfloor\frac{M}{q+1}\rfloor$. 
If $q\ge M\ge 1$, $x^M\Omega_R$ has a limit (zero) as $x$ goes to $0$, then admits a continuous extension.
Otherwise, $ M\ge q+1$ and the differential equation satisfied by $\Omega_R$ gives $$(x^M\Omega_R)' = (Mx^q \text{Id} + R(x)) \cdot x^{M-(q+1)}\Omega_R.$$
From induction hypothesis, $x^{M-(q+1)}\Omega_R$ has a $C^{\lfloor\frac{M}{q+1}\rfloor-1}$ extension.
Then $x^M \Omega_R(x)$ admits a $C^{\lfloor\frac{M}{q+1}\rfloor}$ extension.
\item $\Omega_R$ and $\Omega_R^{-1}$ have the same block diagonal structure than $R$.
\end{enumerate}
\end{proof}

The main result of this section is the following.

\begin{proposition}\label{prop:unlace}\rm
Let $(\xi,\Gamma,(x,\yy))$ be an invariant couple and a system of (TRS) coordinate $(x,\yy): (X,a)\to\mathbb R^{1+n}_0$ of type $(q,N+M,M),$ with:
\begin{enumerate}
\item  $\lfloor \frac{M}{q+1}\rfloor \ge n(q+1+N)+1$;
\item a residual part with spectrum included in $\mathbb R_-^*+i\mathbb R$;
\item a $C^{\infty}$ vestigial part;
\item $(x,\yy)$ has contact of order at least $N+1$ with $\Gamma$.
\end{enumerate} 
Then, 
\begin{enumerate}[(i)]
\item $\xi$ admits a trajectory $(x,\bm\gamma(x)),\; x>0$ which has contact of order at least $N+1$ with $\Gamma$: $\bm\gamma(x)-j_N\Gamma_{\yy}(x) = o(x^N)$.
\item If $(x,\bm\zeta(x)),\; x>0$ is any trajectory of $\xi$ satisfying $\bm\zeta(x)-j_N\Gamma_{\yy}(x)=o(x^N)$, then $\bm\gamma$ and $\bm\zeta$ have flat contact:
 $\forall k\in\mathbb N,\;\bm\gamma(x)-\bm\zeta(x) = o(x^k)$.
  \item For any neighborhood $U$ of $a$, there is an open neighborhood $V\subset U$ of $a$ and a closed connected topological submanifold $S$ of $V\cap\{x>0\}$ such that:
 \begin{enumerate}
 \item $\dim(S)=1+u(D(x))$, where $u(D(x))$ is the unstability index of the exponential part $D(x)$ of the (TRS)-form;
 \item $S$ is locally invariant for $\xi$;
 \item For all $b\in V\cap\{x>0\}$, $b\in S$ if and only if the trajectory of $\xi_{|V}$, parametrized as $(x,\bm \gamma(x))$, $x\in (\alpha,\omega)$, satisfies $\alpha=0$ and 
 $\lim_{x\to 0}\bm\g (x)=0$.
 \end{enumerate}
\end{enumerate} 
\end{proposition}

\begin{proof}
As in Lemma~\ref{prop:key}, the factor $x^e u(x,\yy)$ of the (TRS)-form does not intervene, since the result only involves the 
foliation induced by $\xi$ in the half space $x>0$. So we suppose that
$$\xi=
x^{q+1}\partial_x+\left( \left(D(x)+x^qC \right)\cdot\yy+x^{q+1+N+M}V(x,x^M\yy)\right)\cdot\partial_{\yy},
$$
$D$, $C$, $V$ satisfying the hypothesis.
Up to reorder the coordinate functions $\yy$, we suppose $y_1,\dots,y_m$ carry all dominant rotations. {More precisely}, we suppose that
$$D = \Theta(\text{Diag}(c_1,\dots,c_{m/2}))\oplus\Theta(\text{Diag}(c_{m/2+1},\dots,c_{n_1}))\oplus\text{Diag}(d_{1},\dots,d_{n_2})$$
where:
\begin{enumerate}
\item $c_1,\dots,c_{m/2}\in\mathbb C_{q-1}[x]$, and $\text{ord}_x(\text{Re}(c_j))>\text{ord}_x(\text{Im}(c_j))$ for $j=1,\dots,m/2$;
\item $c_{m/2+1},\dots,c_{n_1}\in\mathbb C_{q-1}[x]$, and $\text{ord}_x(\text{Re}(c_j))\le\text{ord}_x(\text{Im}(c_j))$ for $j=\frac{m}{2}+1,\dots,n_1$;
\item $d_1,\dots,d_{n_2}\in\mathbb R_{q-1}[x]$.
\end{enumerate}
If $m=0$, $(\xi,\Gamma,(x,\yy))$ is in (TRS)-form of type $(q,N,0)$ without dominant rotation and 
{Lemma} \ref{prop:key} applies, which gives the result, taking into account $\Gamma$ has contact at least $N+1$ with the $x$-axis. 

Otherwise, for $l=1,\dots,m/2$, we let 
$$v_l=\text{ord}_x \text{Re} (c_l(x)) \text{ and } b_l(x)=j_{v_l-1}(c_l(x)),$$
so $b_l(x)$ is the initial pure imaginary part of $c_l(x)$. Define the rotational matrix
$R(x) := \Theta(\text{Diag}(b_1(x),\dots, b_{\frac{m}{2}}(x))\oplus \bm 0_{n-m}.$

We apply the transformation $(x,\yy)=U_R(x,\zz)$. The trajectories $(x,\yy(x))$ of $\xi$ included in $\mathbb R_+^*\times\mathbb R^{n}$
are in one-to-one correspondance with the trajectories $(x,\zz(x))$ of {$\eta=U_R^*\xi$} included in $\mathbb R_+^*\times\mathbb R^{n}$, where 
$$ \eta = x^{q+1}\partial_x + \left[(D(x)-R(x) + x^q C )\cdot\zz + x^{q+1+N+M}\Omega_{R}^{-1}(x)\cdot V(x, x^M \Omega_{R}(x) \cdot\zz)\right]\cdot\partial_{\zz}.$$
The crucial point to get this expression is that $C$ and $D(x)$ commute with $\Omega_{R}(x)$ and $\Omega_{R}^{-1}(x)$ 
(Lemma \ref{lem:omega}, (6)).

Let $s\in\{0,\dots,q\}$ be the order of $D(x)-R(x)+x^qC$, and let $$\bm g(x,\zz)=x^M\Omega_{R}^{-1}(x)\cdot V(x, x^M \Omega_{R}(x)\cdot\zz).$$ 
From Lemma \ref{lem:omega} (5),  $\bm g$ admits a $C^{\lfloor\frac{M}{q+1}\rfloor}$ extension on $\mathbb R_+\times\mathbb R^n$. 
The expression of the vector field $\eta/x^s$ is
$$ \eta/x^s = x^{q-s+1}\partial_x + \left[\left(\frac{1}{x^s}(D-R)(x) + x^{q-s} C \right)\cdot \zz + x^{q+1-s+N}\bm g(x, \zz)\right]\cdot \partial_{\zz}.$$
We deduce that $(x,\zz)$ is a system of (TRS)-coordinates of type $(q-s,N,0)$ for $\eta/x^s$, and that:
\begin{itemize}
\item the exponential part $x^{-s}(D-R)(x)$ has no dominant rotation;
\item the residual part $C$ has spectrum included in $\mathbb R^*_-\oplus i\mathbb R$.
\item the vestigial part $\bm g(x,\yy)$  is $C^{\lfloor\frac{M}{q+1}\rfloor}$, and $\lfloor\frac{M}{q+1}\rfloor\ge n(q+1+N)+1$.
\end{itemize} 
Then Lemma \ref{prop:key} applies to $\eta/x^s$, and since $\eta$ and $\eta/x^s$
have the same trajectories, the conclusions of that result applies to $\eta$.
We have that $\eta$ admits a trajectory $(x,\bm \delta(x))$ which has contact of order at least $N+1$ with the $x$-axis, and consequently with $\Gamma$. 
So $(x,\bm \gamma(x)):=(x,\Omega_R(x)\bm\delta(x))$ is a trajectory of $\xi$ which has also contact of order at least $N+1$ with the $x$-axis.
This proves point (i). 
Now, if $(x,\bm \zeta(x))$ is a trajectory of $\xi$ with $\lim\bm\zeta(x)=0$, then  
$(x,\Omega_{R}(x)\bm\zeta(x))$ is a trajectory of $\eta$ that also satisfies $\lim\Omega_{R}(x)\bm\zeta(x))=0$. By Lemma~\ref{prop:key}, (ii),  $\Omega_{R}(x)\bm\zeta(x)$ and $\bm\delta(x)$
have flat contact. Since $\Omega_R$ preserves $\text{ord}_x$, $\bm\zeta(x)$ and $\bm\gamma(x)$ have flat contact, which gives point (ii). Finally, point (iii) is obtained from the corresponding item (iii) of Lemma~\ref{prop:key} applied to $\eta$, taking into account that $\Omega_R$ provides a diffeomorphism from the half-space $\{x>0\}$ to itself, preserving the accumulation of trajectories to the origin. Concerning the dimension of $S$, notice that, according to Definition~\ref{def:no-dominant-rotation}, we have equality of the unstability indices $u(D(x))=u(x^{-s}(D(x)-R(x)))$.
\end{proof}

\section{Final proof}\label{sec:final}

In this section, we prove Theorem~\ref{thm:1}, that is, the existence of trajectories asymptotic to a given formal invariant curve $\G$ of a vector field $\xi$, and Theorem~\ref{thm:2}, which describes the structure (a topological embbeded manifold of positive dimension) of the set of trajectories asymptotic to each half-branch associated to $\G$. The accurate version of these results is Theorem~\ref{th:main-precise} below.  

We recall that a real (irreducible) formal curve $\G$ at {$0\in\R^{m}$} can be defined either with a class of formal parametrizations $\G(t)\in(t\R[[t]])^{m}\setminus\{0\}$ modulo formal change of parameter of the form $t=\alpha(s)\in s\R[[s]]$ with $\alpha'(0)\ne0$, or with a sequence $IT(\Gamma)=(p_k)_{k\ge 0}$ of {\em infinitely near points} or {\em iterated tangents}. This sequence is determined as follows: $p_0=0$, $\G_0=\G$, and, {recursively} for $k>0$, $p_k$ is 
the point in the exceptional divisor of the punctual blowing-up $\pi_{k-1}$ of $p_{k-1}$ where the strict transform $\G_{k}=\pi_{k-1}^*\G_{k-1}$ of $\G_{k-1}$ is centered. We refer to \cite{Wal,Cas} for basics on formal curves and infinitely near points of them.

A formal curve $\G$ has two {\em formal half-branches} $\G^+,\G^-$, defined in the following ways, depending on the chosen definition of $\G$. Considering a parametrization $\G(t)$ of $\G$, the half-branch $\G^\epsilon$, for $\epsilon\in\{+,-\}$, is the equivalence class of $\G^\epsilon(t)=\Gamma(\epsilon t)$ under reparametrization $t=\alpha(s)$ with $\alpha'(0)>0$. In terms of iterated tangents, if we replace the sequence of blowing-ups $\pi_k$ by spherical blowing-ups, $\Gamma$ gives rise to two sequences of {\em oriented iterated tangents} $(q^\epsilon_k)_{k\ge0}$, for $\epsilon\in\{+,-\}$, each one corresponding to a half-branch $\G^\epsilon$ as defined above.
Calculating $IT(\G^\epsilon)$ from $\Gamma^\varepsilon(t)$ involves to know the value of $\lim_{t\to 0} t/|t|$. We adopt
the convention $t>0$.

Although parametrizations and iterated tangents define the same objects, one or the other definition can be more practical to state
a given property. For instance, ``$\G$ is not contained in the singular locus of $\xi$'' 
is concisely stated as $(\widehat{\xi}\circ\G(t))\wedge\G'(t)=0$, where $\Gamma(t)$ is a parametrization of $\Gamma$.
On the contrary, saying that a curve is asymptotic with a half-branch is easily defined with iterated oriented tangents. In general (see \cite{Can-Mou-San1}), a $C^1$ parametrized curve $\g:{(0,c]\to\R^{m}}\setminus\{0\}$ with $\lim_{t\to 0}\g(t)=0$ is said to have {\em (oriented) iterated tangents} if we can associate to $\g$ a sequence of points $IT(\g)=(q_k^+)_{k\ge0}$, where $q^+_0=0$, $\g_0=\g$ and, for $k>0$, $q_k^+=\lim_{t\to 0}\g_k(t)$ where 
$\g_k=\sigma_{k-1}^{-1}\circ\g_{k-1}$ and $\sigma_{k-1}$ is the spherical blowing-up of $q^+_{k-1}$. We say that $\g$ is {\em asymptotic} to the half-branch $\G^\epsilon$ if {$\g$ has iterated tangents and} $IT(\g)=IT(\G^\epsilon)$. This definition clear out the need of a common parameter for $\gamma$ and $\Gamma^{\epsilon}$. When such parameter exists, it can be reformulated, saying the parametrization of $\G^\epsilon$ is an asymptotic expansion of $\gamma$ ``à la Poincaré''. For example, if $\Gamma^{+}=(x,\Gamma^+_{\yy}(x)), x>0$ is a  half-branch of a formal non-singular curve $\G$ at $0\in\R^{1+n}$ and $\gamma = (x,\bm \gamma_{\yy}(x)), x>0$ is a parametrized curve, then $\gamma$ is asymptotic to $\Gamma^+$ if and only if 
\begin{equation}\label{eq:asymptotic}
\forall N\in\mathbb N,\; \|\bm\g_{\yy}(x)-j_N\G^+_{\yy}(x)\| = O(x^{N+1})
\end{equation}
(see for instance \cite[Lemme 4.2]{LeG-Mat-San}). 

Iterated tangents help to define the ``neighborhoods'' of a formal half-branch $\G^+$ where we find asymptotic trajectories. Say the oriented iterated tangents $IT(\G^+)=(q^+_k)_{k\ge0}$ of $\Gamma^+$ are obtained via the sequence of spherical blowing-ups
\begin{equation}\label{eq:equence-real-blow-ups}
\R^{m}\stackrel{\sigma_0}{\leftarrow}M_1\stackrel{\sigma_1}{\leftarrow}M_2\stackrel{\sigma_2}{\leftarrow}\cdots
\end{equation}
A {\em horn neighborhood of $\G^+$} is an open set $V\subset\R^{m}\setminus\{0\}$ 
such that, for some $k\in\mathbb N$, the closure $\overline{V_k}$ of $V_k=(\sigma_0\circ\cdots\sigma_{k-1})^{-1}(V)$ is a neighborhood of $q_k^+$ in $M_k$.
The minimal such $k$ is called \emph{opening} of $V$.
For example, if $\Gamma^+$ has parametrization ${(t,\Gamma_\yy(t))}\in \mathbb R[[t]]^{1+n}$, then for any $\varepsilon>0$, $C>0$, 
\begin{equation}\label{eq:c.h.nbhd}
\Delta(C,\varepsilon)=\{(x,\yy)\in\mathbb R^{1+n};\; 0<x<\varepsilon,\; ||\yy-j_k\Gamma_\yy(x)||<C x^{k}\}
\end{equation} 
is a horn neighborhood of $\Gamma^+$ of opening $k$. Moreover, any 
horn neighborhood of $\Gamma^+$ of opening not smaller than $k$ contains some $\Delta(C,\varepsilon)$.

We can now state our main result in precise terms. 

\begin{theorem}\label{th:main-precise}\rm
Let $\xi$ be a $C^\infty$ vector field in a neighborhood of $a\in\R^{m}$ and let $\G$ be a formal irreducible curve at $a$, invariant for $\xi$ and not contained in the formal singular locus of $\xi$. Let $\G^{\epsilon}$ be a half-branch of $\G$. Then, for any neighborhood $U$ of $a$, {there exists $k_0\in\mathbb N$ such that, for all $k\ge k_0$,}  there is a horn neighborhood $V^\epsilon\subset U\setminus\text{Sing}(\xi)$ of $\G^\epsilon$ of opening $k$, and there is a closed, connected $C^0$-submanifold $S^\epsilon$ of $V^\epsilon$, of positive dimension and locally invariant by $\xi$ such that:
\begin{enumerate}[(i)]
    \item For any $b\in S^\epsilon$, the trajectory $\g$ of $\xi$ through $b$ accumulates to $a$ and is asymptotic to $\G^\epsilon$.
    \item For any $b\in V^\epsilon\setminus S^\epsilon$, the trajectory $\g$ of $\xi$ through $a$ escapes from $V^\epsilon$ in finite time, both positive and negative.
\end{enumerate}
\end{theorem}

Notice that {this result implies the two first statements} of the introduction. Having  positive dimension, $S^{\epsilon}$ is not empty and a trajectory
issued from any point of $S^{\epsilon}$ is asymptotic to $\Gamma$ (so Theorem \ref{thm:1}). Also, any trajectory asymptotic to $\Gamma^{\epsilon}$ will
eventually be included in any horn neighborhood of $\Gamma^{\epsilon}$, in particular $V^{\epsilon}$, so $S^{\epsilon}$ contains the germ of any trajectory asymptotic to $\Gamma^{\epsilon}$ (so Theorem \ref{thm:2}).

We devote the rest of the section to the proof of Theorem~\ref{th:main-precise}. We fix a half-branch, say $\G^+$. Similarly to section~\ref{sec:admissible-for-vf}, the couple $(\xi,\G^+)$ is called an {\em invariant couple} {(although $\G$ might be singular)}. We denote the sequence of spherical blowing-ups giving rise to $IT(\G^+)=\{q^+_k\}_k$ as in (\ref{eq:equence-real-blow-ups}). For each $k$, the composition $\Sigma_k:=\sigma_0\circ\sigma_1\circ\cdots\circ\sigma_k$ provides the {\em transformed} invariant couple $(\xi_k,\G^+_k):=\Sigma_k^*(\xi,\G^+)$, where $\xi_k=\Sigma_k^*\xi$ is the pull-back of $\xi$ by $\Sigma_k$ and $\G^+_k=\Sigma_k^*\G^+$ is given by $IT(\G^+_k)=\{q^+_\ell\}_{\ell\ge k}$. Notice {$\G_k:=\Sigma_k^*\G$} is invariant by $\xi_k$ and not contained in the formal singular locus $\text{Sing}(\widehat{\xi}_k)$. 

\vspace{.2cm}

{\em Step 1. Resolution of $\G$ and adapted coordinates.-} By the definitions of asymptotics and of horn neighborhood so introduced, and since blowing-ups are isomorphisms outside the divisor, it suffices to prove Theorem~\ref{th:main-precise} for the invariant couple $(\xi_k,\G^+_k)$, where $k$ is any arbitrary fixed integer number. 

Using reduction of singularities of formal curves by blowing-ups {(see \cite{Wal})}, we therefore assume that the invariant curve $\G$ is non-singular, that is, we start with an invariant couple $(\xi,\G)$ as considered in section~\ref{sec:admissible-for-vf}. {For notational convenience, we put $m=1+n$ for the dimension of the ambient space and we} adopt all notations and definitions from that and subsequent sections, with the obvious minor modifications needed to handle the half-branch $\G^+$. In particular, a coordinate system $(x,\yy)$ at $a$ is adapted for $\G^+$ if it is adapted for $\G$ ($\Gamma$ is transverse to $\yy=0$) and $\G^+$ is included in the half space $x>0$, that is $\G^+$ has a formal parametrization $(x,\G^+_\yy(x))$, $x>0$.

\vspace{.2cm}

{\em Step 2. Reduction to TRS-form.-} 
Choose adapted coordinates $(x,\yy)$ for $\G^+$. By Theorem~\ref{th:TRS-vf}, (i), there is a composition $\psi=\phi_1\circ\cdots\circ\phi_1$ of admissible coordinate transformations for $(\xi,\G,(x,\yy))$ such that $(\wt{\xi},\wt{\G},(\wt{x},\wt{\yy}))=\psi^*(\xi,\G,(x,\yy))$ is in a TRS-form of type $(q,0,0)$ whose residual part has good spectrum. According to Definition~\ref{def:admissible:coordinates}, the map $\psi$ is written in coordinate as
$$
\psi(\wt x, \wt\yy) = (\wt{x}^\ell,T(\wt{x})\cdot\wt{\yy}+
\bm\alpha(\wt{x})),
$$
where $\ell\ge1$ and the entries of $T$ and $\bm\alpha$ are polynomials in $\wt{x}$.
In particular, $\{x >0\}\subset \psi(\{\wt x>0\})$, so $\Gamma^+$ has a lift $\wt{\G}^+:=\psi^*\G^+$ parametrized as $(\wt x,\wt\Gamma_{\wt\yy}^+(\wt x)),\, \wt{x}>0$, where $(x,\G^+_\yy(x))=\psi(\wt{x},\wt{\G}^+_{\wt\yy}(\wt{x}))$. Also, if $\wt\gamma$ is a trajectory of $\wt{\xi}$ parametrized as  $(\wt{x},\wt{\bm\g}(\wt{x})),\,  \wt{x}>0$, its image $\g=\psi(\wt{\g})$ is a trajectory of $\xi$, parametrized as $$(x,T(x^{1/\ell})\wt{\bm\g}(x^{1/\ell})+
\bm\alpha(x^{1/\ell})), x>0.$$
Considering the reformulation (\ref{eq:asymptotic}) of asymptotic, and the basis (\ref{eq:c.h.nbhd}) of horn neighborhoods, it suffices to prove Theorem~\ref{th:main-precise} for the invariant couple $(\wt{\xi},\wt{\G}^+)$, which is in a TRS-form of type $(q,0,0)$ and whose residual part has good spectrum.

\vspace{.2cm}

{\em Step 3. Existence of asymptotic trajectories.-}
According to the previous step, we consider $(x,\yy)$, a system of TRS coordinates for $(\xi,\G^+)$, of type $(q,0,0)$, with a residual part with good spectrum. 
We can apply Theorem~\ref{th:TRS-vf}, (ii) and reason as in Step 2, then
we can assume that $(\xi,\G)$ is in (TRS)-form of type $(q,M_0,M_0)$, for any given $M_0$.
We choose $M_0$ so that the hypothesis (1)--(4) of  Proposition~\ref{prop:unlace}
are satisfied with $N=0, M=M_0$.
We still denote by $(x,\yy)$ these (TRS) coordinates.
We can now prove the following.
\begin{lemma}\label{lm:one-asymptotic-trajectory}
There exists a trajectory $\g_0= (x,\bm\g_0(x)),\, x>0,$ asymptotic to $\G^+$. 
\end{lemma}
\begin{proof}
Let $(x,\bm\g_0(x)), x>0,$ be a trajectory of $\xi$ such that $\bm\g_0(x)=o(x)$ as provided by Proposition~\ref{prop:unlace}, (i). Let us show that $\g_0$ is asymptotic to $\G^+$.

Given $N\in\N_{\ge1}$, we choose $M=M(N)$ satisfying (1) of Proposition~\ref{prop:unlace}. Then, using Theorem~\ref{th:TRS-vf}, (ii), we take a finite composition of admissible coordinate transformations $\Psi_{N,M}$ such that the transformed invariant couple $(\xi_{N,M},\G_{N,M},(x,\yy_{N,M})) :=\Psi^*_{N,M}(\xi,\G,(x,\yy))$ is in (TRS)-form of type $(q,N+M,M)$ and satisfies hypothesis (1)--(3) of Proposition~\ref{prop:unlace}. We finish with a polynomial translation transformation (but keep {the same} notations) {so that the coordinates $(x,\yy_{N,M})$} also satisfy the hypothesis (4) of that proposition.
Notice that the first coordinate is preserved by $\Psi_{N,M}$, so the system $(x,\yy_{N,M})$ is adapted to the invariant couple $(\xi_{N,M},\G^+_{N,M})=\Psi_{N,M}^*(\xi,\G^+)$.

Applying Proposition~\ref{prop:unlace}, (i), we get, for each $N\ge1$, a trajectory $\g_N=(x,\bm\g_{N}(x))$, $x>0$, of $\xi_{N,M}$ 
which has contact of order  at least $N+1$ with $\G^+_{N,M}$.
Since $\Psi_{N,M}$ is a composition of punctual blowing-ups and affine translations, admissible for $\G$, the image $\zeta_N(x):=\Psi_{N,M}(x,\bm\g_N(x))=(x,\bm\zeta_N(x))$ is a trajectory of $\xi$, tangent to $\G$ of order at least $N+1$.
In particular, since $\lim_{x\to 0}\zeta_N(x)=0$, Proposition~\ref{prop:unlace}, (ii) for $(\xi,\G)$ gives that $\bm\g_0(x)$ and $\bm\zeta_N(x)$ have flat contact, for any $N\ge 1$. We conclude
that $\g_0$ has the same contact order with $\G$ than any $\zeta_N$. Thus, $\g_0$ is asymptotic to $\G^+$.
\end{proof}

{\em Step 4. Submanifold of asymptotic trajectories.-}
To complete the proof of Theorem~\ref{th:main-precise}, we show the existence of the sets $V^+,S^+$ with the stated properties. 
Recall we already assumed that $(\xi,\G^+)$ is in (TRS)-form of type $(q,M_0,M_0)$ in {adapted coordinates $(x,\yy)$}  and that the hypothesis (1)--(4) of Proposition~\ref{prop:unlace} are fulfilled.
 Let $\Psi_{1,M_1}$ be {a} sequence of admissible punctual blowing-ups and polynomial translations 
 such that $(\xi_{1,M_1},\G^+_{1,M_1}, (x,\yy_{1,M_1})):=\Psi_{1,M_1}^*(\xi,\G^+,(x,\yy))$ 
 is in (TRS)-form of type $(q,M_1+1,M_1)$ and also satisfies (1)--(4){; that is, for $N=1$ and $M=M_1$}. 
{Even if it means performing additional punctual blowing ups, $\Psi_{1,M_1}$ can be assumed to contain exactly $k$ of them, for any given $k\ge 1$ larger than a certain $k_0$.}
 
Let $U$ be an initial neighborhood of $a$. From Proposition~\ref{prop:unlace}, (iii) applied  to $(\xi_{1,M_1},\G_{1,M_1})$, we get a neighborhood $V_1$ of the origin of the chart $(x,\yy_{1,M_1})$, contained in $\Psi_{1,M_1}^{-1}(U)$, and a connected, closed $C^0$-submanifold $S_1$ of $V_1\cap\{x>0\}$ such that for any $b\in V_1\cap\{x>0\}$, the trajectory $\delta$ of $\xi_{1,M_1}$ through $b$ satisfies $\alpha(\delta)=0$ if and only if $b\in S_1$. We assume, moreover, that $V_1$ is relatively compact and that $\xi_{1,M_1}(x)>0$ on $\overline{V_1}\cap\{x>0\}$. Thus, taking the parametrization $\delta(x)=(x,\bm\delta(x))$, the trajectory $\delta$ always scapes $V_1$ for positive time, while $\lim_{x\to 0}\bm\delta(x)=0$ if and only if $b\in S_1$.

Define $V^+=\Psi_{1,M_1}(V_1\cap\{x>0\})$ and $S^+=\Psi_{1,M_1}(S_1)$. So $V^+$ is a horn neighborhood of $\G^+$ contained in $U$ and $S^+$ is a closed connected $C^0$-submanifold of $V^+$, locally invariant for $\xi$. Since $V_1$ is relatively compact,  $V^+$ has opening $k$, the number of punctual blowing-up in the composition $\Psi_{1,M_1}$.
Let us check that properties (i) and (ii) in Theorem~\ref{th:main-precise} hold. Take $b\in V^+$ and let {$\g=(x,\bm\g(x))$} be the trajectory of $\xi$ through $b$. Denote by $\wt{\g}=\Psi_{1,M_1}^{-1}(\g)$ the lifted trajectory of $\xi_{1,M_1}$ through $\tilde{b}=\Psi_{1,M_1}^{-1}(b)$. The curve $\g$ escapes $V^+$ for positive time since $\wt{\g}$ escapes $V_1$ for positive time. 

If $b\in S^+$, then $\wt{b}\in S_1$ and hence $\alpha(\wt{\g})=0$, which shows that $\alpha(\g)=0$, too. Moreover, since $k\ge 1$, $\g$ has contact order at least $1$ with $\G$. Applying Proposition~\ref{prop:unlace}, (ii) to $(\xi,\G)$,  $\bm\g(x)$ and $\bm\g_0(x)$ have flat contact, where $\g_0=(x,\bm\g_0(x))$ is the asymptotic trajectory obtained in Step 3. We deduce from (\ref{eq:asymptotic}) that $\g$ is asymptotic to $\G^+$. This gives point (i).

On the contrary, if $b\in V^+\setminus S^+$ then $\alpha(\g)$ cannot be $0$: otherwise, by Proposition~\ref{prop:unlace}, (ii) again, $\g$ would be asymptotic to $\G^+$, and in particular $\wt{\g}$ would accumulate to the origin of the chart $(x,\yy_{1,M_1})$, so, by the properties stated for $V_1$ and $S_1$, we would have $\tilde{b}\in S_1$, which is a contradiction. Then $\g$ escapes $V^+$ for negative time also, since $\overline{V^+}\cap\{x=0\}=\{0\}$, being $V^+$ of opening $k>0$. This gives point (ii).

This finishes the proof of Theorem~\ref{th:main-precise}. \hfill{$\square$}
\section{Non-oscillating trajectories}

In this section, we prove {Theorem~\ref{cor:ana}} which realizes
formal invariant curves of analytic vector fields as subanalytically non-oscillating 
trajectories.
Recall a parametrized curve {$\gamma:(0,c)\to\mathbb R^{m}$}
is subanalytically non-oscillating if its support $|\gamma|=\gamma((0,c))$
intersects any global subanalytic set $X\subset \mathbb R^{m}$ in finitely
many connected components.
When $|\gamma|$ is a half curve at $a\in\mathbb R^{m}$, say
$a=\lim_{t\to 0} \gamma(t)$ and $a\not\in \overline{\gamma(\varepsilon,c)}$ for $\varepsilon >0$,
the notion can be localized
as follows. We say the curve $\gamma$ is subanalytically non-oscillating at $a\in\mathbb R^{m}$
if the intersection of the germ of $|\gamma|$ at $a$ with any subanalytic germ
has a connected representative.
{We prove the following, that is the statement of Theorem~\ref{cor:ana} for formal half-branches}

\vspace{.3cm}
{\bf Theorem \ref{cor:ana}'.}
Let $(\xi,\Gamma)$ be an invariant couple at $a\in\mathbb R^{m}$, suppose $\xi$ is analytic,
and let $\Gamma^{\epsilon}$ be a half-branch of $\Gamma$. Then 
$\xi$ admits a trajectory $\gamma$, asymptotic to $\Gamma^{\epsilon}$
and subanalytically non-oscillating at $a$.

\strut

The proof involves an analytic invariant associated to the formal curve $\Gamma$. In coherence with the terms introduced in section~\ref{sec:TRS}, a (local) blowing-up with center $C$ is {\em admissible} for $\Gamma$ if $C$ does not contain $\Gamma$. This is a sufficient condition to define the lift of $\Gamma$ and, by induction, the notion of an {\em admissible} sequence of blowing-ups.

Given a sequence $\pi= \pi_r\circ\dots\circ\pi_1 : M\to (\mathbb R^{m},a)$ 
of blowing-ups with smooth analytic centers, admissible for $\G$,  
denote by $d_{\pi}(\Gamma)$ the minimal dimension of an analytic 
set $X\subset M$ which contains $\pi^* \Gamma$. We call {\em subanalytic dimension
of $\Gamma$} the minimum of {the} $d_{\pi}(\Gamma)$ when $\pi$ ranges in all
such sequences of blowing-ups, and we denote it by $d(\Gamma)$.
Note that the subanalytic dimension cannot decrease by admissible blowing-ups.

\vspace{.2cm}

\begin{proof}
Choose $\pi,X$ that realize $d(\Gamma)$; that is, 
$\pi= \pi_r\circ\dots\circ\pi_1 : M\to \mathbb R^{m}_p$ is a sequence 
of blowing-ups with analytic smooth centers admissible for $\Gamma$,
and $X\subset M$ is an analytic set of dimension $d(\Gamma)$ which contains {$\pi^*\Gamma$}. 
From Hironaka's Reduction of Singularities for 
analytic sets \cite{Hir} (see also \cite{Bie-M}), there is a composition $\rho: A\to M$ of blowing-ups with
analytic smooth centers such that the strict transform $\rho^*(X)$ {of $X$} is {non-singular}.
Since the centers involved in the resolution have positive codimension
in $X$, the minimality of $d(\Gamma)$ implies that $\rho$ is admissible for $\pi^*\Gamma$.
Let $\Delta=(\pi\circ\rho)^*\Gamma$ {and let $q = \Delta(0)\in A$ be the point where $\Delta$ is centered}. 
The sequence $\pi\circ\rho$ might not be admissible for $\xi$, but 
being only interested by the trajectories of $\xi$ in a neighborhood of $\Gamma$,
we can lift $\xi$ weakly, even if it means to multiply the vector field by an analytic
function whose zero set is the center of each considered blowing-up. We denote by $\zeta\in \text{Der}_{\mathbb R}(\mathcal O(A,q))$
this weak lift. 

We claim $\rho^*X$ is invariant by $\zeta$. Indeed, the tangency locus $T$ of 
$\zeta$ with $X$ is an analytic subset of $X$, and it contains $\Delta$.
By minimality of $d(\Gamma)$, $T$ has dimension larger than $d(\Gamma)$, but since $X$ is
non-singular, it contains no proper analytic subset of dimension $d(\Gamma)$. So $T=X$.

Now $\rho^*X$ is a {regular analytic submanifold}, the restriction $\wt\zeta:=\zeta|_{\rho^*X}$ of $\zeta$ to $X$
is a smooth vector field in $\text{Der}_{\mathbb R}(C^{\infty}(\rho^*X,q))$, and $\Delta$ is a formal curve, invariant by
$\wt\zeta$ that is not included in its formal singular locus. Theorem~\ref{th:main-precise}
applies. There exists a trajectory $\delta$ of $\wt\zeta$ that is asymptotic to the half-branch 
$\Delta^{\epsilon}$ of $\Delta$ which corresponds to $\Gamma^{\epsilon}$. We let $\gamma=(\pi\circ\rho)(\delta)$ be the blow down of $\delta$. Notice that, since $\delta$ is asymptotic to
$\Delta$ and $\Delta$ is not included in the exceptional
divisor of $\pi\circ\rho$, $\gamma$ is truely a curve (and not a point).  

The curve $\gamma$ is by construction a trajectory of $\xi$ asymptotic to $\Gamma^{\varepsilon}$. We claim that $\gamma$ is subanalytically non-oscillating, which will complete the proof.
For this, let $S$ be a subanalytic subset of $\mathbb R^{m}$, and suppose that the germ of $S\cap |\gamma|$ {at $a$} is not empty. We shall prove $|\gamma|\subset S$.
Write $S'=(\pi\circ\rho)^*S\cap \rho^*X$ for the intersection of $\rho^*X$ with the lift of $S$. Since $|\delta|\subset \rho^*X$, the germ of $S'\cap |\delta|$ {at $q$} is not empty. 

We apply Hironaka's Rectilinearization Theorem \cite{Hir}
to the subanalytic set $S'$. There is a covering of a neighborhood of $q$ in $A$ by finitely many semi-analytic sets, and for each of them, a sequence of blowing-ups which 
transforms $S'$ into an analytic set. Call $U$ the 
semi-analytic set of this partition which contains $\Delta$, and $\sigma:B\to U$ the corresponding sequence of blowups. Again, the minimality of $d(\Gamma)$ implies that
$U$ contains an open subset of $\rho^*X$ and $\sigma$ is admissible for $\Delta$. The lift $\sigma^*S'$ is an analytic set which contains infinitely many points of {$\sigma^{-1}(\delta)$}. Then it contains $\sigma^*\Delta$. Since $d(\sigma^*\Delta)\ge d(\Gamma) = \dim X$ {(where $\sigma^*\Delta$ is considered as a formal curve in the ambient $m$-manifold $B$), we get that} $\sigma^*S'\cap (\rho\circ\sigma)^*X$ has dimension not smaller than {$\dim X=\dim (\rho\circ\sigma)^*X$}. But $(\rho\circ\sigma)^*X$ being {non-singular}, it has no proper analytic subset of the same dimension. Then $\sigma^*S'=(\rho\circ\sigma)^*X$, and pushing forward, $\rho^*X\subset S'$, then $|\delta|\subset S'$ and therefore $|\gamma|\subset S$.
\end{proof}



\end{document}